\documentclass[11pt,a4paper,reqno]{amsart}
\usepackage{graphicx}
\usepackage{amsmath} 
\usepackage{amssymb}
\usepackage{mathtools}
\usepackage{color}
\usepackage{graphicx}
\usepackage{booktabs}
\usepackage{color}
\usepackage{braket}
\usepackage{algorithm}
\usepackage{cases}
\usepackage{array}
\usepackage[
  colorlinks=true
]{hyperref}

\hypersetup{
	linkcolor={red!50!black},
	citecolor={blue!50!black},
	urlcolor={blue!80!black}
}
\usepackage[labelformat=simple]{subcaption}

\usepackage[shortlabels]{enumitem}
\usepackage[foot]{amsaddr}
\usepackage{amsopn}
\usepackage[nameinlink]{cleveref}
\crefformat{equation}{\textup{#2(#1)#3}}

\usepackage[utf8]{inputenc}
\usepackage{amsmath,amssymb,amsfonts}
\usepackage{latexsym}
\usepackage{mathrsfs}
\usepackage{bbm}
\usepackage{xcolor}
\usepackage{url}
\usepackage{graphicx}
\usepackage{stmaryrd}
\usepackage{enumitem}
\SetSymbolFont{stmry}{bold}{U}{stmry}{m}{n}
\usepackage{graphicx}

\usepackage{etoolbox}  
\usepackage{refcount}  

\makeatletter
\newif\if@intheorem
\@intheoremfalse
\pretocmd{\theorem}{\@intheoremtrue}{}{}
\apptocmd{\endtheorem}{\@intheoremfalse}{}{}
\makeatother 

\newcommand{\storeEquationInTheorem}[1]{%
  \expandafter\gdef\csname orig@eq@#1\endcsname{\theequation}%
}

\newcommand{\origEquation}[2][]{%
\begin{equation}
\if\relax\detokenize{#1}\relax
  \tag*{\csname orig@eq@#2\endcsname} 
\else
  \tag*{#1} 
\fi
\label{#2}
\end{equation}
}

\providecommand{\restatetheorem}{}  
\renewcommand{\restatetheorem}[2][]{%
  \begingroup
  \renewcommand{\thetheorem}{\getrefnumber{#2}}%
  \begin{theorem}[#1]
}

\providecommand{\restatetheoremend}{}  
\renewcommand{\restatetheoremend}{
  \end{theorem}
  \endgroup
}

\usepackage{dsfont}

\usepackage{subcaption}
\usepackage{booktabs}

\usepackage{algpseudocode}
\usepackage{algorithm}

\usepackage[shortcuts]{extdash}

\usepackage{enumitem}
\setlist[enumerate]{leftmargin=.5in}
\setlist[itemize]{leftmargin=.5in}

\newcommand{\kernel}{\text{\sf{k}}}

\crefname{theorem}{Theorem}{Theorems}
\crefname{lemma}{Lemma}{Lemmas}
\crefname{proposition}{Proposition}{Propositions}
\crefname{corollary}{Corollary}{Corollaries}
\crefname{assumption}{Assumption}{Assumptions}
\crefname{claim}{Claim}{Claims}

\crefname{definition}{Definition}{Definitions}
\crefname{remark}{Remark}{Remarks}
\crefname{example}{Example}{Examples}
\crefname{idea}{Idea}{Ideas}
\crefname{hypothesis}{Hypothesis}{Hypotheses}

\newcommand{\one}{\mathbbm{1}}

\newcommand{\<}{\langle}
\renewcommand{\>}{\rangle}

\newcommand{\wh}{\widehat}

\newcommand{\R}{\mathbb R}    
\newcommand{\N}{\mathbb N}    

\newcommand{\calB}{\mathcal B}         
         
\newcommand{\calD}{\mathcal D}

\newcommand{\calI}{\mathcal I}         
         
\newcommand{\calK}{\mathcal K}         
         
\newcommand{\calM}{\mathcal M}         
\newcommand{\calN}{\mathcal N}

\newcommand{\calV}{\mathcal V}         
         
\newcommand{\calX}{\mathcal X}         
\newcommand{\calY}{\mathcal Y}

         \newcommand{\bE}{\mathbb E}

         \newcommand{\bH}{\mathbb H}

         \newcommand{\bN}{\mathbb N}
         
         \newcommand{\bP}{\mathbb P}
         
         \newcommand{\bR}{\mathbb R}

\newcommand{\veps}{\varepsilon}

\newcommand{\vphi}{\varphi}

\newcommand{\bbK}{\mathbf K}
\newcommand{\bbk}{\mathbf k}

\allowdisplaybreaks

\newtheorem{theorem}{Theorem}[section]
\newtheorem*{theorem*}{Theorem}
\newtheorem{remark}[theorem]{Remark}
\newtheorem{proposition}[theorem]{Proposition} %
\newtheorem{corollary}[theorem]{Corollary}  %
\newtheorem{lemma}[theorem]{Lemma} %

\newtheorem{assumption}[theorem]{Assumption}

\title[Koopman for stochastic dynamics: error bounds for kernel EDMD]{
	Koopman for stochastic dynamics: error bounds for kernel extended dynamic mode decomposition
}

\author{Maximiliano Hertel$^{1}$}
\address{$^1$Optimization-based Control Group, Institute of Mathematics,
Technische Universit\"at Ilmenau, Germany
(\textsc{\{maximiliano.hertel, friedrich.philipp, karl.worthmann\}@tu-ilmenau.de}).}%
\author{Friedrich M.~Philipp$^{1}$}
\author{Manuel Schaller$^{2}$}
\address{Chemnitz University of Technology, Germany
(\textsc{manuel.schaller@math.tu-chemnitz.de}).}%
\author{Karl Worthmann$^{1}$}

\thanks{M. Hertel gratefully acknowledges funding by Carl-Zeiss-Stiftung, project KI-MSO-O: AI-supported analysis, modeling and synthesis (design) of organoids, P2024-02-019. 
F.M. Philipp is grateful for the  support of the German Research Foundation (DFG), project numbers 554600805 and 519323897.
K. Worthmann gratefully acknowledges funding by the DFG, project ALeSCo: Active Learning for Systems and Control, project number 535860958..}

\begin{document}
\begin{abstract}
    We prove $L^\infty$-error bounds for kernel extended dynamic mode decomposition (kEDMD) approximants of the Koopman operator for stochastic dynamical systems. To this end, we establish Koopman invariance of suitably chosen reproducing kernel Hilbert spaces and provide an in-depth analysis of the pointwise error in terms of the data points. The latter is split into two parts by showing that kEDMD for stochastic systems involves a kernel regression step leading to a deterministic error in the fill distance as well as Monte Carlo sampling to approximate unknown expected values yielding a probabilistic error in terms of the number of samples. We illustrate the derived bounds by means of Langevin-type stochastic differential equations involving a nonlinear double-well potential.
\end{abstract}

\maketitle


\smallskip
\noindent \textbf{Mathematics subject classications.}      37M99, %
    47B32, %
    65C05, %
    65D12 %

\bigskip

\section{Introduction}\label{ch:intro}
Nowadays, the Koopman operator is an established mathematical and computational tool for modeling, analysis, and control of dynamical systems, see, e.g., the collection~\cite{MaurSusu20} or the recent overview papers~\cite{brunton:budisic:kaiser:kutz:2022,strasser2026overview}. 
The underlying idea is to consider nonlinear dynamics through the lens of observable functions using the linear, but infinite-dimensional Koopman operator, see the original work by B.O.~Koopman~\cite{Koop31,KoopNeum32} and the seminal paper~\cite{Mezi2005} by I.~Mezi{\'c}. This enabled, among others, to leverage modern techniques from spectral analysis to investigate the complex dynamical behaviour of the original nonlinear system~\cite{rowley2009spectral}.

Recently, data-driven approaches like dynamic mode decomposition (EDMD; \cite{WillKevr15,Colb23}) were proposed to approximate the Koopman operator, see, e.g., the recent overview~\cite{Colb24}, which has lead to a plethora of successful applications, including fluid flows~\cite{bagheri:2013, mezic:2013}, climate forecasting~\cite{hogg:fonoberova:mezic:2020,azencot:erichson:lin:mahoney:2020}, molecular dynamics~\cite{nuske:klus:2023}, power systems technology~\cite{susuki:mezic:raak:hikihara:2016}, quantum systems~\cite{klus:nuske:peitz:2022}, robotics~\cite{mamakoukas:castano:tan:murphey:2019,bruder:fu:gillespie:remy:vasudevan:2020,haggerty:banks:kamenar:cao:curtis:mazic:hawkes:2023,rosenfelder:bold:eschmann:eberhard:worthmann:ebel:2024}, and human-centered applications such as grip force prediction for robotic rehabilitation~\cite{bazina:kamenar:fonoberova:mezic:2025}.
In EDMD, snapshots of the underlying system are employed to approximate the Koopman operator. It has been shown in~\cite{KordaMezi18} that EDMD converges in the infinite-data and infinite-dictionary (infinitely many observable functions) limit, corresponding to a complexity solvability index of two~\cite{colbrook2022foundations}. Then, first finite-data results on the estimation error for ergodic sampling for deterministic systems were provided in \cite{Mezi22} followed by~\cite{ZhanZuaz23} analyzing the $L^2$-projection error onto a finite element dictionary, cf.\ also the works \cite{NuskPeit23,PhilScha24} providing error bounds for stochastic and control systems. 

More recently, kernel EDMD approximation schemes, which allow a data-informed choice of the dictionary, were proposed, see, e.g., \cite{WillRowl15,klus:nuske:hamzi:2020}.  
Several works have already addressed $L^2$\=/error bounds for kernel EDMD: \cite{PhilScha24:kernel} provides these in terms of the Mercer eigenvalues of the chosen kernel in the case of ergodic sampling along a trajectory, while the bounds in \cite{PhilSchaWort2025,OlguOsseRami2025}, are derived for i.i.d.\ sampling. The work \cite{KohnPhil24} shows first pointwise (i.e., $L^\infty$) error bounds for approximations using Wendland kernels. 
However, while kernel EDMD is successfully used for stochastic systems, pointwise error bounds  are still missing --~to the best of the authors' knowledge.

In this paper, we propose a readily applicable Monte Carlo framework for estimation of kernel matrices to compute approximations of the Koopman operator for stochastic systems. Then, we provide an in-depth error analysis of the scheme.
To this end, we establish Koopman invariance of (fractional) Sobolev spaces based on novel sufficient conditions on a Radon-Nikodym density associated with the stochastic dynamical system. Then we split the approximation error into a deterministic and a probabilistic part, where the latter results from Monte Carlo estimation of the expected value. Combining these ingredients, we then present novel full approximation error bounds for kernel EDMD for stochastic systems using Matérn and Wendland kernels.
Finally, we present numerical simulations for systems governed by SDEs as a proof of concept.

The paper is organized as follows. In \Cref{ch:setting}, we introduce the Koopman operator for stochastic systems, propose a kernel EDMD approach for its approximation and state the main result in \Cref{thm:total_error_bound}, providing a pointwise total error bound for the kernel EDMD approximant with high probability. To prove this error bound, we first derive necessary and sufficient conditions for invariance of Sobolev spaces under the stochastic Koopman operator in \Cref{ch:hs_invar}, cf. \Cref{thm:hs_frac_charac}. Then we prove the approximation error bounds in \Cref{ch:err_bounds} and illustrate them using numerical examples in \Cref{ch:numerics}. \Cref{ch:conclusion} concludes and gives an outlook. \\[.3cm]

\noindent \textbf{Notation.} For a set $\Omega \subset \R^n$, we denote the Borel sigma algebra on $\Omega$ (w.r.t.\ the standard topology) by $\calB_\Omega$. The Lebesgue measure of a Borel set $\Delta \in \calB_\Omega$ is denoted by $|\Delta|$; $d x$ or $d y$ indicate integrals w.r.t.\ Lebesgue measure.  
$\chi_A$ denotes the characteristic/indicator function of a set $A \subset \R^n$. The Radon-Nikodym density of a measure $\mu$ w.r.t.\ the Lebesgue measure will be denoted by $\frac{d \mu}{d x}$. 
The identity operator on some space is denoted by $\calI$, the identity matrix especially is denoted by $I$. For Banach spaces $\calX$ and $\calY$, $L(\calX, \calY)$ denotes the space of bounded linear operators from $\calX$ to $\calY$. The norm of $T\in L(\calX,\calY)$ is denoted by $\|T\|_{\calX\to\calY}$. We write $L^p(\Omega)$, $p \in [1, \infty]$, for the usual Lebesgue spaces with respect to Lebesgue measure, $C_b(\Omega)$ for the space of bounded and continuous functions, and $C_c^\infty(\Omega)$ for the space of infinitely differentiable and compactly supported functions, respectively, on $\Omega$; occasionally, we just write $L^p$, $C_b$, or $C_c^\infty$ if the domain is clear from context.
For~$f:\Omega \to \R$, we denote $\|f\|_\infty := \sup_{x\in \Omega} |f(x)|$.  
The $\alpha$-th weak derivative of a function $f: \Omega \to \bR$, $\alpha \in \bN^n$, is denoted by $D^\alpha f$. In the case of a function $f: \Omega^L \to \bR$, $L \in \bN$, taking multiple arguments, we write $D_l^\alpha f$ for the $\alpha$-th derivative of $f$ ``w.r.t.\ the $l$-th variable'', $l \in \{1, \dots, L\}$. 

\section{Koopman operator for stochastic systems: kernel EDMD and error analysis}\label{ch:setting}

Let $\Omega\subset\R^n$ be a bounded domain with Lipschitz boundary and 
$\calB_\Omega$ be the Borel sigma algebra on~$\Omega$. 
Moreover, let $X$ and $Y$ be $\Omega$-valued random variables. 
It is well known (see, e.g., \cite{Klen2008}) that there exists a function $\rho : \Omega\times\calB_\Omega\to [0,1]$ such that
\begin{enumerate}[label=(\roman*)]
	\item $\rho_x := \rho(x,\cdot)$ is a probability measure for each $x\in\Omega$;
	\item For each $\Delta\in\calB_\Omega$, the function $\rho(\cdot,\Delta): \Omega\to [0,1]$, is measurable; 
	\item For each $\Delta \in \calB_\Omega$: $\rho(x,\Delta) = \bP(Y\in\Delta\,|\,X=x)$ for~$\bP^X$-a.e.\ $x\in\Omega$.
\end{enumerate}
A function with the first two properties is called a \emph{(Markov) transition kernel} or \emph{probability kernel}; the third property makes $\rho$ a \emph{regular version} of the conditional distribution of $Y$ given $X$ (see, e.g., \cite{Kall2001}). In what follows, we consider stochastic dynamical systems of the form 
\begin{align}\label{e:stoch_sys}
	Y \mid \{X = x\} \sim \rho_x,\quad x\in \Omega,
\end{align}
where the distribution of $Y$ is conditioned on the initial value $X = x$. 
Specifically, one might consider a discrete-time Markov process, possibly given by the flow of a stochastic differential equation (SDE) for a time step $T > 0$, i.e., 
\begin{align}\label{e:sde_chap2}
	d X_t &= b(X_t) \,dt\, + \sigma(X_t) \,d W_t
\end{align}
with drift vector field $b: \R^n \to \R^n$, diffusion vector field $\sigma: \R^n \to \R^{n \times n}$, and $n$\=/dimensional Brownian motion $W_t$.

\subsection{The Koopman operator for stochastic systems}

The \emph{Koopman operator} $\calK$ associated with $X$ and $Y$ maps functions $f$ on $\Omega$ to functions $\calK f$ on $\Omega$ and is defined as
\begin{align}\label{eq:koop_identity}
	(\calK f)(x) := \bE\big[f(Y)\,|\,X=x] = \int_\Omega f\,d\rho_x,\quad x\in\Omega.
\end{align}
In the next section, we elaborate necessary and sufficient conditions for the invariance of Sobolev spaces $H^s(\Omega)$ under the Koopman operator. These spaces are of special interest since they arise as native spaces, or \emph{reproducing kernel Hilbert spaces} (RKHS), of certain kernels for which there exists an extensive approximation framework, see~\Cref{subch:kedmd}, that allows for error bounds (\Cref{ch:err_bounds}). 
Before analyzing the Koopman operator on Sobolev spaces, we require the invariance of the space of square-integrable functions $L^2(\Omega) = H^0(\Omega)$ under $\calK$.
\begin{lemma}\label{l:l2_invariance}
	If the condition
	\begin{equation} 
		\label{e:L1} 
		\exists C>0 : \quad\int\rho_x(\Delta)\,dx\,\le\,C|\Delta| \qquad\,\forall\Delta\in \calB_\Omega, 
	\end{equation}
	holds, then $\calK L^2(\Omega)\subset L^2(\Omega)$, and $\calK : L^2(\Omega)\to L^2(\Omega)$ is a bounded linear operator.
\end{lemma} 
\begin{proof}
	Indeed, Cauchy-Schwarz yields $|(\calK f)(x)|^2\le (\calK f^2)(x)$, so that for simple functions $f = \sum_{j=1}^mc_j\one_{\Delta_j}$ with mutually disjoint sets $\Delta_j$,
	\[
	\int |(\calK f)(x)|^2\,dx \le \int (\calK f^2)(x)\,dx = \sum_{j=1}^mc_j^2\int\rho_x(\Delta_j)\,dx\,\le\,C\sum_{j=1}^mc_j^2|\Delta_j| = C\|f\|_{L^2(\Omega)}^2.
	\]
	By density of the simple functions in~$L^2(\Omega)$, it follows that $\calK$ is well defined on $L^2(\Omega)$ and that $\|\calK f\|_{L^2(\Omega)}\le C^{1/2}\|f\|_{L^2(\Omega)}$ for all $f\in L^2(\Omega)$.
\end{proof}
Note that, by the Fubini-Tonelli theorem, \eqref{e:L1} implies $L^1$-invariance of the Koopman operator. More generally, it can be shown that for a bounded domain~$\Omega$, $L^1$-invariance (w.r.t.\ Lebesgue measure) implies $L^q$-invariance for~$q > 1$. 

From now on, we consider the Koopman operator acting on functions $\psi \in L^2(\Omega)$.

\begin{remark}
	While we consider $L^2(\Omega) = L^2_{{\pmb{\lambda}}}(\Omega)$ with the Lebesgue measure~$\pmb{\lambda}$, our results can be extended to other reference measures --~depending on the underlying stochastic system. In the case of $Y = X_{T}$ and $X = X_0$ for some Markov process $(X_t)_{t \geq 0}$, the invariant measure $\mu$ of the process would be an intuitive choice, see, e.g., \cite{NuskPeit23}. 
\end{remark}

\subsection{Recap on deterministic systems: kernel EDMD and error bound}\label{subch:kedmd}

Let $\Omega \subset \bR^n$ as before. In the deterministic case, the dynamics are given by
\begin{align*}
	x^+ = F(x), 
\end{align*}
with $F: \Omega \to \Omega$, which may be considered as a particular case of the stochastic case \eqref{e:stoch_sys}. In particular, we may consider the SDE~\eqref{e:sde_chap2} with $\sigma \equiv 0$ giving rise to the ordinary differential equation $\dot x = b(x)$. 
In this case, the transition kernel is given by $\rho_x (A) = \delta_{F(x)}(A)$, such that the Koopman identity~\eqref{eq:koop_identity} reduces to 
\begin{align*}
	\calK \psi(x) = \psi(F(x)),
\end{align*}
i.e., $\calK \psi = \psi \circ F$. As in \cite{KohnPhil24}, we assume that $F : \Omega\to\Omega$ is a diffeomorphism satisfying $\inf_{x\in\Omega}|\det F'(x)| = c > 0$ and $\sup_{x\in\Omega}\|F'(x)\| < \infty$ . While the latter condition corresponds to the Lipschitz continuity of the drift function in \eqref{e:sde_chap2}, the assumption on the Jacobi determinant is a special case of \eqref{e:L1}. Applying the change of variables formula gives
\begin{align*}
	\int \delta_{F(x)}(\Delta) \, dx 
	= \left| F^{-1}(\Delta) \right| = \int_\Delta |\det (F^{-1})'(y)| \, dy \leq \frac{1}{c}\, |\Delta|,
\end{align*}
i.e., \eqref{e:L1} also holds in the deterministic case.

We briefly recall \emph{Extended Dynamic Mode Decomposition} (EDMD; \cite{WillKevr15}) and its extension \emph{kernel EDMD}, which are well-established numerical methods to approximate the Koopman operator $\calK$ from data, see also~\cite{KutzBrun16} or the recent overview article~\cite{Colb24}. 
EDMD requires two ingredients: 
First, let $d \in \bN$ \textit{snapshots} $(x_i, y_i)$ with $y_i = F(x_i)$, $i \in \{1,\dots,d\}$, be given. Then, we partition the snapshot data into the sets
\begin{align*}
	\calX = \{x_1, \dots, x_d \} \quad\text{ and }\quad
	\calY = \{y_1, \dots, y_d\}.
\end{align*}
Second, let $M \in \bN$ observable $L^2$-functions $\psi_1, \dots, \psi_M$ be given, for which we use the notation $\Psi(x) = (\psi_1(x), \dots, \psi_M(x))^\top$, and define the data matrices
\begin{align*}
	\Psi_\calX = (\Psi(x_1), \dots, \Psi(x_d)) \in \bR^{M \times d}, \qquad 
	\Psi_\calY = (\Psi(y_1), \dots, \Psi(y_d)) \in \bR^{M \times d}.
\end{align*}
The subspace $\calV = \mathrm{span}\{\psi_1, \dots, \psi_M\} \subset L^2(\Omega)$ spanned by the observables is finite\=/dimensional.
From the Koopman identity, we get 
\begin{align}\label{eq:koop_psi_map}
	\calK \psi_i(x_j) = \psi_i(F(x_j)) = \psi_i(y_j) \qquad\forall\, (i,j) \in \{1,\ldots,M\}\times\{1, \dots, d\},
\end{align}
that is, the Koopman operator maps the components of $\Psi_\calX$ to the components of $\Psi_\calY$. Since the Koopman operator acts \emph{linearly} on observables, we seek a data-driven linear approximation on~$\calV$ given by a matrix $K \in \R^{M \times M}$ \emph{approximately} satisfying~\eqref{eq:koop_psi_map}, that is,
\begin{align}\label{eq:matrix_equation_koop_edmd}
	\Psi_\calY \approx  K \Psi_\calX.
\end{align}
Assuming that $\Psi_\calX$ has full row rank (which is equivalent to strong linear independence of the observables $\psi_1, \dots, \psi_M$ w.r.t.\ Lebesgue measure, see~\cite[Lemma C.3]{PhilScha24}), an approximation can be computed as a solution to the least squares problem
\begin{align}\label{eq:koop_regression_problem}
	\min_{K \in \bR^{M \times M}} \| K \Psi_\calX - \Psi_\calY \|_F^2.
\end{align}
Direct calculations using the normal equations yield the EDMD approximation
\begin{align}\label{eq:compression_approx_edmd_v2}
	K = \Psi_\calY \Psi_\calX^\dagger = \Psi_\calY \Psi_\calX^\top (\Psi_\calX \Psi_\calX^\top)^{-1}.
\end{align}
As illustrated above, EDMD crucially depends on two ingredients: a finite number of data pairs $\{(x_i, y_i) \mid y_i = F(x_i),\,i \in \{1, \dots, d\} \}$ and the choice of the observables. 
Convergence in the infinite-data limit has been shown in~\cite{KordaMezi18} and further specified w.r.t.\ the estimation error (number of data pairs~$d$) in~\cite{Mezi22,NuskPeit23,PhilScha24}. 
The projection error is typically analyzed using finite-element observables as proposed in~\cite{ZhanZuaz23} and extended to control systems in~\cite{SchaWort23}. 
However, the latter only yields $L^2$-error bounds. 

A particular variant of EDMD using a data-informed set of observables is \emph{kernel EDMD} (kEDMD; \cite{WillRowl15}), sharing two crucial features: First, from an application point of view, the choice of the observable functions is performed in a data-driven manner, eliminating fiddly and time-consuming tuning. Second, the rich mathematical structure of reproducing kernel Hilbert spaces allows for \emph{pointwise deterministic} error bounds. 

For~$\Omega \subset \bR^n$, consider a strictly positive definite symmetric kernel $k: \Omega \times \Omega \to \bR$ that generates a unique Hilbert space $\left(\bH, \langle\cdot, \cdot\rangle_\bH\right)$ of bounded and continuous functions on $\Omega$ for which it is a \emph{reproducing kernel} \cite{Aron1950}.\footnote{We note that kernels are usually defined globally on $\bR^n \times \bR^n$, as is the case for Wendland and Matérn kernels we consider later.} In particular, for all $\psi \in \bH$ the \emph{reproducing property}
\begin{align}\label{eq:rkhs_prop}
	\psi(x) = \langle \psi, \Phi_x\rangle_\bH \quad \forall x \in \bR^n
\end{align}
holds, where $\Phi_x\in\bH$, defined by $\Phi_x(y) = \kernel(x,y)$, is commonly termed the \emph{canonical feature} of $\kernel$ at $x \in \Omega$. In kEDMD, these features play a central role:
We choose or sample points $\calX = \{ x_1, \dots, x_d\}$ as before and let $y_i = F(x_i)$. The dictionary is given by the canonical features $\Phi_x = \kernel(x, \cdot)$, $x \in \Omega$, at the data points, i.e.
\begin{align*}
	\calV = \calV_\calX := \mathrm{span}(\calD) = \mathrm{span}\{ \Phi_{x_1}, \dots, \Phi_{x_d} \}.
\end{align*}
Clearly $M = d$ holds, i.e., we have exactly as many observables~$\Phi_{x_i}$ as data points $x_i$. The data matrix entries are easily calculated as $(\Psi_\calX)_{i,j} = \psi_i(x_j) = \Phi_{x_i}(x_j) = \kernel(x_i, x_j)$ and $(\Psi_\calY)_{i,j} = \psi_i(y_j) = \Phi_{x_i}(F(x_j)) = \kernel(x_i, F(x_j))$. We write them as \emph{kernel matrices} 
\begin{align*}
	\mathbf{K}_\calX := \mathbf{K}_{\calX, \calX} = (\kernel(x_i, x_j))_{i,j=1}^d = \Psi_\calX, \qquad 
	\mathbf{K}_{\calX, F(\calX)} := (\kernel(x_i, F(x_j))_{i,j=1}^d = \Psi_\calY.
\end{align*}
Plugging the kernel matrices into \eqref{eq:koop_regression_problem} yields the regression problem
\begin{align} \label{eq:kedmd_regression_problem}
	\min_{K \in \bR^{d \times d}} \left\| K \bbK_\calX - \bbK_{\calX, F(\calX)} \right\|_F^2.
\end{align}
Assuming that the data points are pairwise distinct, $\mathbf{K}_\calX$ is symmetric positive definite (hence invertible) due to the symmetry and strict positive definiteness of the underlying kernel $\mathbf{k}$.\footnote{We explicitly speak of \emph{strict} positive definiteness above, since positive definiteness of kernels is often defined via positive \emph{semi-definiteness} of kernel matrices in the literature, see~\cite{KohnPhil24}.} The regression problem \eqref{eq:kedmd_regression_problem} now reduces to (exactly) solving a linear system, whose solution $K = \mathbf{K}_{\calX, F(\calX)} \mathbf{K}_\calX^{-1}$ corresponds to \eqref{eq:compression_approx_edmd_v2}. This matrix not only represents the action of the Koopman operator for our data, but gives the matrix representation of $S_\mathcal{X} \mathcal{K}\vert_{\calV}$ in the basis of the canonical features, where $S_\calX := S_{V_\calX}$ denotes the orthogonal projection onto $\calV$. Specifically, for some $g \in \bH$, the projection onto $\calV_\calX$ is given by 
\begin{align}\label{eq:projection_formula}
	S_\calX g = \sum_{i=1}^d \left( \bbK_\calX^{-1} g_\calX \right)_i \Phi_{x_i} = g_\calX^\top \bbK_\calX^{-1} \bbk_\calX
\end{align} 
with $g_\calX := (g(x_1), \dots, g(x_d))^\top$ and $\bbk_\calX := (\Phi_{x_1}, \dots, \Phi_{x_d})^\top$, see~\cite[Proposition 2.1]{KohnPhil24}. For~$f \in \calV_\calX$, there exists $\alpha_f \in \bR^d$ such that $f = \alpha_f^\top \bbk_\calX$. With $g = \calK|_{\calV_\calX} f$ we get $g_\calX = \bbK_{\calX, F(\calX)}^\top \alpha_f$ and 
\begin{align}\label{eq:kedmd_approx_alpha}
	S_\calX \calK|_{\calV_\calX} f = S_\calX g = g_\calX^\top \bbK_\calX^{-1} \bbk_\calX = \alpha_f^{\top} \bbK_{\calX, F(\calX)} \bbK_\calX^{-1} \bbk_\calX.
\end{align}
This means the \emph{kEDMD}-approximant $K = \mathbf{K}_{\calX, F(\calX)} \mathbf{K}_\calX^{-1}$ \cite[Proposition 3.2]{KohnPhil24} maps (the coefficients of) $f \in \calV_\calX$ to (the coefficients of) $S_\calX \calK|_{\calV_\calX} f$ linearly (in the respective canonical bases). To extend this operator to the entire space $\bH$, one may simply add another orthogonal projection from the right, i.e., replace $f \in \calV_\calX$ by $S_\calX g$ for any $g \in \bH$. From \eqref{eq:projection_formula}, we get the representation $S_\calX g = \alpha_{S_\calX g}^\top \bbk_\calX = g_\calX^\top \bbK_\calX^{-1} \bbk_\calX$. In total, this yields the Koopman approximant 
\begin{align}\label{eq:koop_kedmd_approximant}
	\wh \calK = S_\calX \calK S_\calX: \bH \mapsto V_\calX,
\end{align}
which, using \eqref{eq:kedmd_approx_alpha}, may be evaluated by 
\begin{align}\label{eq:koop_edmd_projection_complete}
	\wh \calK f = S_\calX \calK S_\calX f = S_\calX \calK|_{\calV_\calX} (S_\calX f) = \alpha_{S_\calX f}^\top \bbK_{\calX, F(\calX)} \bbK_\calX^{-1} \bbk_\calX = f_\calX^\top \bbK_\calX^{-1} \bbK_{\calX, F(\calX)} \bbK_\calX^{-1} \bbk_\calX.
\end{align}
The main result regarding this deterministic Koopman approximation \cite[Theorem 3.4]{KohnPhil24} may be summarized as follows.
\begin{theorem}
	\label{thm:main_prior_work}
	Suppose that $\calK\bH\subset\bH$. Then $\calK : \bH\to\bH$ is a bounded operator, and the approximant $\wh\calK = S_\calX \calK S_\calX$ satisfies  
	\begin{align*}
		\sup_{\|f\|_\bH\leq 1}\| (\calK - \widehat{\calK})f\|_{\infty} \leq\, (\|\calK\|_{C_b\to C_b}+\|\calK\|_{\bH\to \bH})\cdot \|\calI - S_\calX \|_{\bH \to C_b}. 
	\end{align*}
\end{theorem}
The above error bound depends on two crucial features:
\begin{enumerate}
	\item \textbf{Invariance property}. First, it is assumed that the Koopman operator is indeed a bounded operator on $\bH$. Due to closedness of the Koopman operator on RKHS, it suffices to show that $\calK \bH\subset \bH$ (see~\cite[Lemma 3.1]{KohnPhil24} and \Cref{prop:RKHS_inv} below).
	\item \textbf{Approximation error}. 
	Second, for kEDMD, the regression problem is indeed solved by interpolation and the interpolation error $\|\calI - S_\calX \|_{\bH \to C_b}$ has to be controlled; this is usually done in terms of the \emph{fill distance}~$h_\calX$ of the data set $\calX$, i.e., one usually has
	\begin{align}\label{eq:errorbound_fill}
		\|\calI - S_\calX\|_{\bH \to C_b} \leq C h_\calX^{r} \qquad\text{with}\qquad h_\calX = \underset{z \in \Omega \; x \in \calX}{\sup \min}\,\| z - x \|_2
	\end{align}
	for a suitable rate $r>0$ and $C\geq 0$ using existing results from kernel-based approximation theory, see e.g.\ \cite[Section 11.2]{Wend04}.
\end{enumerate}
The invariance property (1) is a subtle detail that has to be verified rigorously and cannot be taken for granted. In particular, there are RKHS for which $\mathcal{K}\in L(\bH,\bH)$ implies that the underlying dynamics is affine linear, such as the Gaussian RKHS, see~\cite{GonzAbud23}. For (fractional) Sobolev spaces $H^s$ (see~\Cref{subch:hs_invar_sobol_spaces} for an introduction), however, which are RKHS if $s > n/2$ and are generated by Wendland or Matérn kernels (see below for a precise definition), Koopman invariance was proven in \cite[Theorem 4.2]{KohnPhil24} in the deterministic setting. The central assumption is that the flow is smooth enough, such that the Koopman composition operator preserves Sobolev regularity. On the other hand, invariance of a RKHS under the Koopman operator is sufficient for boundedness, as shown in \cite[Lemma 3.1]{KohnPhil24} for the deterministic case, and analogously in \Cref{prop:RKHS_inv} below for the stochastic case. 

\begin{proposition}\label{prop:RKHS_inv}
	Let $\bH\subset C_b(\Omega)$ be a RKHS that is invariant under $\calK$. Then $\calK : \bH\to\bH$ is bounded.
\end{proposition}
\begin{proof}
	Since $\bH\subset C_b(\Omega)$, the corresponding kernel is bounded. Let $f_n,f,g\in\bH$ such that $f_n\to f$ and $\calK f_n\to g$ in~$\bH$. Recall the canonical features $\Phi_x$ at~$x \in \Omega$. Then, the reproducing property \eqref{eq:rkhs_prop} yields
	\begin{align*}
		g(x) &= \<g,\Phi_x\>_\bH = \lim_{n\to\infty}\<\calK f_n,\Phi_x\>_\bH = \lim_{n\to\infty}\calK f_n(x) = \lim_{n\to\infty}\int_\Omega f_n(y)\,d\rho_x(y).
	\end{align*}
	Now, $f_n(y) = \<f_n,\Phi_y\>_\bH\to\<f,\Phi_y\>_\bH = f(y)$ for all $y\in\Omega$ as $n\to\infty$ by continuity of the scalar product, and $|\<f_n, \Phi_y \>_\bH|\le\|f_n\|_\bH\|\Phi_y\|_\bH\le C$ for all $n\in\N$. Thus, by dominated convergence, $g(x) = \calK f(x)$, which proves that $\calK$ is closed and, thus, bounded.
\end{proof}

To sum up, the assumption of \Cref{thm:main_prior_work} is met for Sobolev RKHS, and bounds on the projection error in the form \eqref{eq:errorbound_fill} are well-known for many generating kernels, leading to a \emph{pointwise} and \emph{deterministic} error bound for deterministic systems.

\subsection{Kernel EDMD for stochastic systems}\label{subch:stoch_kedmd}
In constrast, the stochastic Koopman operator propagates conditional expectations of observables. Here, the counterpart of the approximate identity \eqref{eq:matrix_equation_koop_edmd} takes the form
\begin{align}\label{eq:matrix_equation_koop_edmd_stoch}
	K \Psi_\calX\approx \Psi_\calY = \big(\bE[\Psi_i(Y) \mid X = x_j]\big)_{i,j=1}^{M ,d}.
\end{align}
Consequently, for the finite-dimensional approximation $K \in \bR^{M \times M}$ in \eqref{eq:matrix_equation_koop_edmd_stoch}, setting
\begin{align*}
	Z^+ := \Psi(Y) \approx K Z, \qquad Z := \Psi(X),
\end{align*}
defines a \emph{linear surrogate system in feature space}.
The surrogate system describes the evolution of the \emph{conditional mean} of the lifted observables, i.e., $\bE[Z^+\,|\, X] \approx K Z = K \Psi(X)$.
Here, probability enters through the conditional expectation: the surrogate system approximates the expected evolution of the lifted observables, while stochastic fluctuations around this mean remain. However, the conditional expectations in \eqref{eq:matrix_equation_koop_edmd_stoch} are usually not readily available in applications with unknown dynamics, such that the analogue to the regression problem \eqref{eq:koop_regression_problem} cannot be directly solved. Instead, the data matrix $\Psi_\calY$ has to be approximated from observations of the stochastic system. For this, we build upon the rich theory of Monte Carlo methods to perform an entry-wise approximation. We draw $m_j \in \bN$ i.i.d.\ samples $y_j^{(l)} \sim  \rho_{x_j}$, $l = 1, \dots, m_j$, and invoke the law of large numbers such that
\begin{align}\label{eq:kxy_mc_lln_approx}
	\bE[\Psi_i(Y) \mid X = x_j] \approx \frac{1}{m_j} \sum_{l =1}^{m_j} \Psi_i(y_j^{(l)}).
\end{align}
For the kernel case, we then have, in view of $\Psi_i = \kernel(x_i,\cdot)$,
\begin{align}\label{eq:matrix_equation_koop_edmd_stoch_MC}
	K \mathbf{K}_\calX \approx \mathbf{K}_{\calX, \calY}^\mathrm{MC} \qquad \mathrm{with}\qquad (\mathbf{K}_{\calX, \calY}^\mathrm{MC})_{ij}
	= \frac{1}{m_{j}} 
	\sum_{l=1}^{m_{j}}
	\kernel(x_i, y_j^{(l)}) 
\end{align}
This yields the matrix approximant $K = \bbK_{\calX, \calY}^{\mathrm{MC}} \bbK_\calX^{-1}$, which corresponds to the operator approximant 
\begin{align}\label{eq:koop_edmd_projection_stochastic}
	\widehat{\mathcal{K}}^\mathrm{MC} f = f_{\mathcal{X}}^\top\,
	\bbK_\calX^{-1}\, \mathbf{K}_{\calX,\calY}^\mathrm{MC}\, 
	\bbK_\calX^{-1}\,
	\bbk_{\mathcal{X}},
\end{align}
mirroring the deterministic approximant from \eqref{eq:koop_edmd_projection_complete}. However, the entry-wise Monte Carlo approximation $\bbK_{\calX,\calY}^\mathrm{MC}$ introduces a further estimation error since \eqref{eq:kxy_mc_lln_approx} only holds \emph{approximately}. In the ideal case where $\Psi_\calY = \bbK_{\calX, \calY}$ is known exactly, this estimation error vanishes, and the error bound from \Cref{thm:main_prior_work} holds. 

\begin{remark}  
	Another possible design choice more in line with the covariance operator approach circumvents the estimation of the expected value in \eqref{eq:kxy_mc_lln_approx} under certain assumptions on the underlying stochastic system and sampling scheme (see~\cite{NuskPeit23,PhilSchaWort2025} for details). In the setting of a governing stochastic differential equation, the authors either sample the design points in~$\calX$ from the invariant measure $\mu$ or assume ergodicity and sample from a single long trajectory $x_0, x_1, \dots, x_d$.
	In this paper, we pursue the Monte Carlo approach illustrated above, which requires weaker assumptions on the system but relies on i.i.d.\ samples.
\end{remark}

\subsection{Main contribution: error analysis of kEDMD for stochastic systems}\label{subch:main_contrib}

In this part, we summarize our main results. In addition, we provide a roadmap for the subsequent sections, where the derivations are presented in full detail. Our first main result is an abstract error bound for general RKHS resulting in a \emph{pointwise} and \emph{probabilistic} error bound on the approximation with kEDMD for stochastic systems. 
\begin{theorem}[$L^\infty$-error bound]\label{thm:total_error_bound}
	Let $\bH$ be an RKHS of functions over $\Omega$, and assume $\calK\in L(\mathbb{H}, \mathbb{H})$ and $\calK\in L(C_b(\Omega),C_b(\Omega))$. Furthermore, for~$\mathcal{X} = \{x_1,\ldots,x_d\} \subset \Omega$, let $m_j = m \in \bN$ i.i.d.\ samples from every $\rho_{x_j}$, $j \in \{1, \dots, d\}$, be given.
	Then, for every $\veps > 0$ satisfying $m \geq \frac{2\|\kernel\|_\infty}{\veps^2}\log(2d)$ and  for all $f \in \bH$, $\| f \|_\bH \leq 1$, we have the pointwise error bound
	\begin{equation}\label{eq:total_err}
		\big\| \big(\calK - \widehat{\calK}^{\mathrm{MC}}\big) f\big\|_{\infty} \leq \big[ 1 + \|\calK\|_{\bH\to \bH} \big] \|\calI - S_\calX \|_{\bH \to C_b} + \veps \cdot \|\kernel\|_\infty^{1/2} \sqrt{\max_{v\in\one}v^\top \mathbf{K}_\calX^{-1}v} 
	\end{equation}
	with probability at least $1 - 2d \exp\left(-\frac{m\veps^2}{2\|\kernel\|_\infty}\right)$, where~$\one = \{v\in\R^d : |v_i|=1\}$.
\end{theorem}

The proof of the above theorem (cf.\ \Cref{subch:proof_total_approx}) consists of two parts. First, the finite number of data points and observables given by $\calX \subset \Omega$ introduces a \emph{deterministic error}. Second, we can only estimate the entries of the ``propagation matrix'' $\bbK_{\calX, \calY}^\mathrm{MC}$ which are conditional expected values. Our estimator $\bbK_{\calX, \calY}^\mathrm{MC}$ is random itself, resulting in a \emph{probabilistic} estimation error. Thus, we split the approximation error as follows:
\begin{align*}
	\|\calK - \wh\calK^{\mathrm{MC}}\|_{\bH \to C_b}
	\leq \underbrace{\|\calK - \wh\calK\|_{\bH \to C_b}}_{\text{deterministic error}} + \underbrace{\|\wh\calK - \wh\calK^{\mathrm{MC}}\|_{\bH \to C_b}}_{\text{probabilistic error}}.
\end{align*}

In \Cref{thm:total_error_bound}, the parameter $m$ is chosen sufficiently large (depending on $\varepsilon$ and $d$) to ensure that the stated lower bound on the probability is nonnegative. More precisely, for fixed $d$ and $\varepsilon > 0$, choosing $m \ge \frac{2\|\kernel\|_\infty}{\varepsilon^2}\log(2d)$ guarantees that the probability bound is meaningful.
In practice, the parameters $d$ and $m$ play different roles: the number of design points $d$ controls the deterministic approximation error, while $m$ determines the accuracy of the Monte Carlo estimation. For a prescribed accuracy level $\varepsilon$, one therefore chooses $d$ large enough to control the deterministic error and then selects $m = m(d,\varepsilon)$ sufficiently large to achieve the desired error bound with high probability.
We refer the reader to \Cref{subch:numerics_ou} and \Cref{app:numerical_considerations} for further discussion on the interplay between $d$, $\varepsilon$, and $m$.

The assumption $\calK \in L(C_b(\Omega),C_b(\Omega))$ is a regularity condition on the transition kernel $\rho_x$, which is characterized and discussed in \Cref{subch:det_approx}, and is independent of the continuity of the RKHS kernel $\kernel$ defining $\bH$. Besides, the two main components in \Cref{thm:total_error_bound} are the operator norms (which are assumed to be finite) and the projection error governed by $\calI - S_\calX$. These two components strongly rely on the chosen kernel and thus deserve separate treatments which constitute the second and third main result of this work.

\medskip
\noindent \textbf{Boundedness of the Koopman operator.} 
Invariance of reproducing kernel Hilbert spaces under Koopman operators is a non-trivial issue. Hence, in \Cref{ch:hs_invar} we provide characterizations and sufficient conditions for the invariance of Sobolev spaces. Here, we present the most important and practically relevant result, which follows from the more general  result of \Cref{thm:hs_frac_charac}.

\begin{theorem}\label{thm:diff_density_hs_invariance}
	If the transition kernel $\rho_x$ is given by a density $p: \Omega \times \Omega \to [0, \infty)$, i.e., $\rho_x(\Delta) = \int_\Delta p(x,y) \,dy$  for all $x \in \Omega$ and $\Delta \in \calB_\Omega$,
	that satisfies
	\begin{align}\label{e:Hs_condition_density_chap2}
		\sum_{|\alpha|\leq \ell} \left\Vert D_1^\alpha p(\cdot,\cdot) \right\Vert_{L^2(\Omega \times \Omega)} \leq C
	\end{align}
	for some constant $C > 0$ and $\ell\in \N$, then $\mathcal{K}\in L(H^s(\Omega),H^s(\Omega))$ for all $s \in (0, \ell]$. 
\end{theorem}
We point out that Koopman invariance, as shown in~\cite{KohnPhil24} for suitably chosen reproducing kernel Hilbert spaces, was extended to stochastic systems in the preprint~\cite{OlguOsseRami2025}, under the stronger regularity assumption $p \in C^m(\Omega \times \Omega)$ for the transition density. 
On the one hand, we generalize this result by significantly weakening the assumptions. 
On the other hand, we also provide necessary conditions and, thus, a complete characterization --~even without assuming that a density for the transition kernel $\rho_x$ exists.

\medskip
\noindent \textbf{Projection error bounds in terms of the fill distance.} In this work, we apply our results to two popular choices of kernels, \emph{Wendland kernels} \cite{Wend04} and \emph{Matérn kernels} \cite{BerlThom2004} which are in particular consistent with the above mentioned invariance as they generate \emph{Sobolev spaces} (see~\cite{Adams2003} and \Cref{subch:hs_invar_sobol_spaces} for an introduction). We briefly present the error bounds for Wendland kernels, as the bounds for Matérn kernels (see~\Cref{app:matern_kernels} for an introduction) are structurally very similar. 

The Wendland kernel $\kernel_{n,l}: \bR^n \times \bR^n \to \bR$ is bounded, continuous and strictly positive-definite (see~\cite{Wend04} and \Cref{app:wendland_kernels} for details). For~$\sigma_{n,l} = \frac{n+1}{2} + l$, it generates the (fractional) Sobolev space $H^{\sigma_{n,l}}(\bR^n)$ as RKHS with equivalence of the standard Sobolev norm and the induced (native) RKHS norm. Restricted to a bounded domain~$\Omega \subset \bR^n$ with Lipschitz boundary, as is of interest in the setting of this paper, the Wendland kernel generates the corresponding Sobolev space $H^{\sigma_{n,l}}(\Omega)$ \cite[Theorem 4.1]{KohnPhil24}. Results from kernel approximation theory allow for a more explicit bound on the projection error in terms of the fill distance.

\begin{theorem}[$L^\infty$-error bound for Wendland native spaces] \label{thm:total_err_wendland}
	Let $l \in \bN$, define $\sigma_{n,l} = \frac{n+1}{2}+l$ and denote by $\bH$ the Sobolev space $H^{\sigma_{n,l}}(\Omega)$ induced by the Wendland kernel $\kernel_{n,l}$ as RKHS. Assume that $\calK \in L(\bH, \bH)$ and $\calK \in L(C_b(\Omega), C_b(\Omega))$. Furthermore, for~$\mathcal{X} = \{x_1,\ldots,x_d\} \subset \Omega$, let $m_j = m \in \bN$ i.i.d.\ samples from every $\rho_{x_j}$, $j \in \{1, \dots, d\}$, be given.
	Then, for every $\veps > 0$ satisfying $m \geq \frac{2 C_2^2}{\veps^2}\log(2d)$ and all $f \in \bH$ with $\|f\|_\bH \leq 1$, we have the pointwise error bound
	\begin{equation}\label{eq:total_err_wendland}
		\begin{aligned}
			\left\| (\calK - \calK^{\mathrm{MC}})f\right\|_{\infty} \leq\, &C_1\,(1+\|\calK\|_{\bH\to \bH})\cdot h_\calX^{l + 1/2} + \veps \cdot C_2\cdot\sqrt{\max_{v\in\one}v^\top \mathbf{K}_\calX^{-1}v}
		\end{aligned}
	\end{equation}
	with probability at least $1 - 2d\exp\left(-\frac{m\veps^2}{2 C_2^2}\right)$, where $C_1, C_2 > 0$ are constants only depending on the kernel $\kernel_{n,l}$, and $\one = \{v\in\R^d : |v_i|=1\}$.
\end{theorem}

Note that the bounds of \Cref{thm:total_error_bound,thm:total_err_wendland} depend on the smallest eigenvalue of the kernel matrix~$\bbK_\calX$. 
Tikhonov regularization is a means to control the latter, for the price of introducing a deterministic bias. \Cref{ch:err_bounds} addresses the impact of regularization on the deterministic (\Cref{prop:err_tikh}), probabilistic (\Cref{cor:mc_error1_tikh}), and total error bounds (\Cref{cor:reg_total_err_bound}).

\section{\texorpdfstring{$H^s$}{Hs}-invariance under the Koopman operator}\label{ch:hs_invar}

Since Wendland (and Matérn) kernels generate the Sobolev spaces $H^s(\Omega)$ as their respective native spaces, we study their invariance under the stochastic Koopman operator. As those spaces are defined via partial integration, a condition for preserving weak derivatives under Koopman operator action is paramount. The central ingredient for this is the following signed measure associated with the transition kernel $\rho_x$. For~$g \in L^\infty(\Omega)$, define $\nu_g: \calB_\Omega \to \bR$ by  
\begin{align}\label{eq:nu_g}
	\nu_g(\Delta) := \int_\Omega \rho_x(\Delta)g(x)\,dx,\qquad \Delta \in \calB_\Omega.   
\end{align}
The following lemma links the signed measure $\nu_g$ to the action of the dual of the Koopman operator on $L^1(\Omega)$ and $L^2(\Omega)$. Its proof is given in \Cref{subch:hs_invar_thms}. 

\begin{lemma}\label{l:dual_nu_dens}
	Assume that \eqref{e:L1} holds and consider $\calK_1 : L^1(\Omega)\to L^1(\Omega)$ as well as $\calK_2 : L^2(\Omega)\to L^2(\Omega)$. Then, the measure $\nu_g$ defined in \eqref{eq:nu_g} is absolutely continuous w.r.t.\ the Lebesgue measure, and its density $\frac{d\nu_g}{dx}\in L^\infty(\Omega)$ satisfies 
	\[
	\calK_1^*g = \calK_2^*g = \frac{d\nu_g}{dx}\quad\mathrm{for \ all}\quad g\in L^\infty(\Omega),
	\]
	where $\calK_1^\ast$ and $\calK_2^\ast$ denote the adjoints w.r.t.\ the respective dual pairing. 
\end{lemma}

In \Cref{subch:hs_invar_thms}, we state our main result \Cref{thm:hs_frac_charac} regarding Sobolev space invariance, which turns out to be closely linked to the properties of such densities $\frac{d\nu_g}{dx}$. By~\Cref{prop:RKHS_inv}, this RKHS invariance implies boundedness of the Koopman operator, which proves crucial for the error bound provided later in \Cref{ch:err_bounds}. \Cref{subch:hs_invar_sobol_spaces} briefly introduces relevant notation and defines integer-order and fractional-order Sobolev spaces. The concluding \Cref{subch:hs_invar_proof} is dedicated to the proofs of the results of \Cref{subch:hs_invar_thms}. 

\subsection{Conditions for \texorpdfstring{$H^s$}{Hs}-invariance and boundedness}\label{subch:hs_invar_thms}
The following theorem provides conditions for~$H^s$-invariance in different settings. 
While (i) fully characterizes $H^\ell$-invariance for integer order $\ell \in \bN$, (ii) gives a sufficient condition for~$H^s$-invariance with real-valued $s \geq 0$. (iii) covers the important special case of a (weakly) differentiable transition density and was presented as \Cref{thm:diff_density_hs_invariance} in \Cref{ch:setting}. The proof is deferred to \Cref{subch:hs_invar_proof}.
\begin{theorem}[Koopman invariance of Sobolev spaces]\label{thm:hs_frac_charac} Assume that \eqref{e:L1} holds.
	\begin{itemize}
		\item[{\rm (i)}] Let $\ell \in \bN_{\geq 1}$. Then $H^\ell(\Omega)$ is invariant under $\calK$ if and only if there is $C\geq 0$ such that
		\begin{align}\label{e:hl_charac}
			\sum_{|\alpha|\leq \ell} \left\| \frac{d\nu_{D^\alpha \vphi}}{dx} \right\|_{H^{-|\alpha|}} \leq C \| \vphi \|_{L^2} \qquad \forall \vphi \in C_c^\infty.
		\end{align}
		\item[{\rm (ii)}] If \eqref{e:hl_charac} holds, then $H^s(\Omega)$ is invariant under $\calK$ for all $s \in (0, \ell]$.
		\item [{\rm (iii)}] Assume the transition kernel is given by a density $p: \Omega \times \Omega \to [0, \infty)$  such that 
		\begin{align*}
			\rho_x(\Delta) = \int_\Delta p(x,y) \,dy \quad \mathrm{for\ all}\quad x \in \Omega, \,\Delta \in \calB_\Omega.
		\end{align*}
		If for some constant $C > 0$ and $\ell\in \bN_{\geq 1}$, we have
		\begin{align}\label{e:Hs_condition_density}
			\sum_{|\alpha|\leq \ell} \left\Vert D_1^\alpha \,p(\cdot,\cdot) \right\Vert_{L^2(\Omega \times \Omega)} \leq C,
		\end{align}
		then $H^s(\Omega)$ is invariant under $\calK$ for all $s \in (0, \ell]$. 
	\end{itemize}
\end{theorem}

Note that part (iii) of \Cref{thm:hs_frac_charac} is precisely the statement of \Cref{thm:diff_density_hs_invariance}. In the following result, we exemplarily verify its sufficient condition \eqref{e:Hs_condition_density} for~$H^s$-invariance in case of a particular smooth transition density. Similar results can be obtained for higher-order Sobolev spaces or other (truncated) distributions with smooth densities on bounded domains~$\Omega$. 

\begin{proposition}[Truncated Normal Distribution] \label{prop:trunc_norm}
	For~$n = 1$ and $\Omega = (0,1) \subset \bR$, assume $Y\mid \{X = x\} \sim \rho_x := \mathcal{N}_{(0,1)}(x, 1)$, where $\mathcal{N}_{(0,1)}(x, 1)$ denotes the (doubly) truncated normal distribution on $(0,1)$ with location parameter $x \in (0,1)$ and scale parameter set to $1$. Then, $H^1(\Omega)$ is invariant under the associated Koopman operator that acts on $f \in H^1(\Omega)$ via $f \in H^1\bigl((0,1)\bigr)$ via $\calK f (x) = \bE[f(Y) \mid X = x]$.
\end{proposition}
\begin{proof}
	The doubly truncated normal distribution $\rho_x = \mathcal{N}_{(0,1)}(x, 1)$ is defined as the distribution of a random variable $Z \sim \calN(x,1)$ conditioned on the event $0 < Z < 1$. To apply part (iii) of \Cref{thm:hs_frac_charac}, we study the density of the transition kernel~$\rho_x$. Denote the cumulative distribution function of the standard normal distribution by $F(\xi) = \frac{1}{\sqrt{2\pi}} \int_{-\infty}^\xi e^{-t^2/2} \,dt$. Then, the distribution of $Y$ conditioned on $\{X = x \}$ with $x \in (0,1)$ has density $p_x = p(x, \cdot): (0,1) \to [0, \infty)$ given by 
	\begin{align*}
		p(x,y) = \frac{\frac{1}{\sqrt{2\pi}} e^{-(y-x)^2/2}}{F(1-x) - F(-x)} = \frac{G(x,y)}{H(x)}
	\end{align*}
	The numerator and denominator both depend smoothly on $x$. With $F'(\xi) = \frac{1}{\sqrt{2\pi}}e^{-\xi^2/2}$, we compute the derivative with respect to $x$:
	\begin{align*}
		\partial_1 p(x,y) = \frac{\partial_1 G(x,y) H(x) - G(x,y) H'(x)}{H(x)^2} = G(x,y) \left( (y-x) - \frac{H'(x)}{H(x)^2}\right).
	\end{align*}
	Note that $0 < C_1 \leq H(x) \leq C_2$ since $H(x) = \bP(U \in (-x, 1 - x))$ for some $U \sim \calN(0,1)$ and $x \in (0,1)$. With $G(x,y) \leq \frac{1}{\sqrt{2\pi}}$ for all $x,y \in (0,1)$ and $\vert H'(x) \vert \leq \frac{1}{\sqrt{2\pi}}$, we obtain the bound
	\begin{align*}
		\vert \partial_1 p(x,y) \vert \leq \frac{1}{\sqrt{2\pi}} \left(1 + \frac{1}{\sqrt{2\pi} C_1}\right)
	\end{align*}
	for all $x, y \in (0,1)$. Thus, condition \eqref{e:Hs_condition_density} is satisfied with $C = \frac{1}{\sqrt{2\pi}} \left(1 + \frac{1}{\sqrt{2\pi} C_1}\right)$. By part (iii) of \Cref{thm:hs_frac_charac}, $H^1(\Omega)$ is invariant under the associated Koopman operator $\calK$. 
\end{proof}

\subsection{Sobolev Spaces}\label{subch:hs_invar_sobol_spaces}

For self-containedness, we briefly recapitulate the notion of Sobolev spaces following \cite{Adams2003}. Readers who are familiar with Sobolev spaces may only skim this section for notation or skip it entirely. 

Consider the bounded domain~$\Omega \subset \R^n$ with Lipschitz boundary and the space $C_c^\infty(\Omega)$ of infinitely differentiable and compactly supported functions on $\Omega$.
We recall the notion of the weak derivative of a function: Consider a locally integrable function $f \in L_\mathrm{loc}^1(\Omega)$; that is, $f$ is integrable on every compact set $K \subset \Omega$, for which $f \in L^2(\Omega)$ is clearly sufficient. If there exists a function $g \in L_{\mathrm{loc}}^1(\Omega)$ satisfying 
\begin{align}\label{eq:sobolev_weak_deriv_def}
	\int_\Omega f(x) \, (D^\alpha \vphi)(x)\,dx = (-1)^{|\alpha|} \int_\Omega g(x) \,\varphi(x)\,dx
\end{align}
for every test function $\varphi \in C_c^\infty(\Omega)$, we call $g$ the $\alpha$-th weak derivative of $f$ and write $g = D^\alpha f$. Using the notion of a weak derivative, given some $\ell \in \bN$ we introduce
\begin{align*}
	H^\ell(\Omega) = W^{\ell, 2}(\Omega) = \{ f \in L^2(\Omega): D^\alpha f \in L^2(\Omega)~\text{for all}~ \alpha \in \bN^n~\text{with}~|\alpha| \leq \ell\}
\end{align*}
as the space of $L^2$-functions whose weak derivatives up to degree $\ell$ exist and lie in~$L^2$. If the domain~$\Omega$ is clear from context, we briefly write $H^\ell = H^\ell(\Omega)$. Equipped with the norm
\begin{align*}
	\| f \|_{H^\ell} = \left( \sum_{|\alpha| \leq \ell} \| D^\alpha f \|_{L^2}^2 \right)^{1/2},
\end{align*}
$H^\ell$ is a Banach space.
Its dual space with respect to the $L^2$-scalar product will be denoted by $H^{-\ell}(\Omega)$ or $H^{-\ell}$. We briefly remind the reader that we have the following chain of inclusions with continuous embeddings $H^1\,\subset\,C_b\,\subset\,L^2\,\subset\,C_b^* = \calM\,\subset\,H^{-1}$, where $C_b^* = \calM$ is the space of all finite regular signed measures on $\Omega$. 

Sobolev spaces can be generalized to real-valued regularity order $s \geq 0$. Those \emph{fractional Sobolev spaces} arise, for example, as reproducing kernel Hilbert spaces of certain kernels.
There are two predominant definitions of fractional Sobolev spaces in the literature, either via the Fourier transform or via weak derivatives and interpolation between integer-order Sobolev spaces. For a domain with Lipschitz boundary, as is $\Omega \subset \bR^n$, those definitions coincide by~\cite[Theorem 3.30]{McLe2000}, so we proceed with the second definition which is more suitable for the remainder of this paper. 
Given $s \geq 0$, write $s = \lfloor s \rfloor  + r$ with $r \in [0,1)$ and define the norm
\begin{align}\label{e:frac_sobolev_norm}
	\|f \|_{H^s(\Omega)}^2 := \sum_{| \alpha | \leq\, s} \| D^\alpha f \|_{L^2(\Omega)}^2 + \chi_{(0,1)}(r) \int_\Omega \int_\Omega \frac{\left| D^\alpha f(x) - D^\alpha f(y) \right|^2}{\| x - y \|_2^{n + 2r}}\, dx \,dy.
\end{align}
Under this definition, the Sobolev space $H^s(\Omega) = \{ f \in L^2(\Omega): \| f \|_{H^s(\Omega)} < \infty\}$ is an interpolation space between $H^{\lfloor s\rfloor}(\Omega)$ and $H^{\lceil s\rceil}(\Omega)$. More precisely, $H^s(\Omega) = [H^{\lfloor s\rfloor}(\Omega), H^{\lceil s\rceil}(\Omega)]_r$.

\subsection{Proof of \texorpdfstring{\Cref{l:dual_nu_dens} and \Cref{thm:hs_frac_charac}}{Lemma 3.1 and Theorem 3.2}}\label{subch:hs_invar_proof}

We start by proving the important auxiliary result from the beginning of this section. While \Cref{l:dual_nu_dens} essentially follows from \cite[Appendix B]{PhilScha24}, we provide a concise proof below for self-containedness.
\begin{proof}[Proof of \Cref{l:dual_nu_dens}]
	Let $g \in L^\infty(\Omega)$ and $\Delta \in \calB_\Omega$. If \eqref{e:L1} holds, the definition of \eqref{eq:nu_g} immediately gives
	\[
	|\nu_g(\Delta)|\,\le\,\|g\|_{L^\infty}\int_\Omega \rho_x(\Delta)\,dx\,\le\,C\|g\|_{L^\infty}|\Delta|.
	\]
	Hence, $\nu_g$ is absolutely continuous w.r.t.\ Lebesgue measure with density $\frac{d\nu_g}{dx}\in L^\infty(\Omega)$. Moreover, for~$f\in L^1(\Omega)$ and the corresponding Koopman operator $\calK_1: L^1(\Omega) \to L^1(\Omega)$, we have
	\[
	\<\calK_1f,g\>_{L^1,L^\infty} = \int_\Omega \calK_1f\cdot g\,dx = \int_\Omega \int_\Omega f\,d\rho_x g(x)\,dx = \int_\Omega f\,d\nu_g = \left\<f,\frac{d\nu_g}{dx}\right\>_{L^1,L^\infty}.
	\]
	This shows $\calK_1^*g = \frac{d\nu_g}{dx}$. Finally, since $|\Omega|<\infty$ implies the inclusion $L^\infty(\Omega) \subset L^2(\Omega)$, we have for~$f,g \in L^\infty(\Omega)$:
	\[
	\<\calK_2^* f,g\>_{L^1,L^\infty} = \<\calK_2^* f,g\>_{L^2} = \<f,\calK_2 g\>_{L^2} = \<f,\calK_1 g\>_{L^1,L^\infty} = \<\calK_1^* f,g\>_{L^1,L^\infty}.
	\]
	It follows that $\calK_2^*f = \calK_1^*f$.    
\end{proof}

Building on \Cref{l:dual_nu_dens}, the remainder of the proof of \Cref{thm:hs_frac_charac} consists of five steps. 
\begin{proof}[Proof of \Cref{thm:hs_frac_charac}] 
	
	\noindent \textbf{Step 1. $H^1$-invariance.}
	As the first step, we prove the statement of  for~$\ell=s=1$.
	
	We start with the sufficiency of invariance in (i): Assume that $H^1$ is invariant under $\calK$ and fix $\alpha \in \bN^n$ with $|\alpha| = 1$. Then $\calK : H^1\to H^1$ is bounded by~\Cref{prop:RKHS_inv}. Specifically, there exists $C_1 > 0$ such that 
	\begin{align}\label{eq:h1_bound}
		\| \calK f \|_{H^1} \leq C_1 \| f \|_{H^1} 
	\end{align} 
	for all $f \in H^1$. Moreover, \Cref{l:dual_nu_dens} gives $\frac{d \nu_{D^\alpha \vphi}}{dx} = \calK_2^\ast (D^\alpha \vphi)$, so we study the $H^{-1}$-norm of the latter. For all $f\in H^1$ and all $\vphi\in C_c^\infty$, we have
	\begin{equation}\label{eq:h1_bound_density_ineq}
		\begin{aligned}
			|\<f,\calK_2^* (D^\alpha \vphi) \>_{H^1,H^{-1}}|
			&= \left|\<\calK f, D^\alpha\vphi\>\right| = \left|\< D^\alpha(\calK f),\vphi\>\right|\,\le\,\|D^\alpha(\calK f)\|_{L^2}\|\vphi\|_{L^2}\\
			&\le\,\|\calK f\|_{H^1}\|\vphi\|_{L^2}\,\le\,C_1\|f\|_{H^1}\|\vphi\|_{L^2}.
		\end{aligned}
	\end{equation}
	The first equality follows from the definition of the dual operator, while the second equality leverages the definition of the weak derivative $D^\alpha$. The inequalities follow from the Cauchy-Schwarz inequality, the definition of the norm in the Sobolev space $H^1$ and \eqref{eq:h1_bound}. The inequality above implies
	\begin{align*}
		\left\|\frac{d \nu_{D^\alpha \vphi}}{dx}\right\|_{H^{-1}} = \|\calK_2^* (D^\alpha \vphi)\|_{H^{-1}} 
		&= \sup\left\{\left|\<f,\calK_2^* (D^\alpha \vphi) \>_{H^1,H^{-1}}\right| : f\in H^1,\,\|f\|_{H^1}=1\right\}\,\\&\le\,C_1\|\vphi\|_{L^2}  
	\end{align*}
	for every $|\alpha| = 1$. For~$|\alpha| = 0$, we have $H^0 = L^2$ and can directly apply the Cauchy-Schwarz inequality and \Cref{l:l2_invariance} to get $\left|\<\calK f, \vphi\>\right|  \leq \| \calK f\|_{L^2} \|\vphi\|_{L^2} \leq C_2 \|f\|_{L^2} \|\vphi\|_{L^2}$ and hence
	\begin{align*}
		\left\|\frac{d \nu_{\vphi}}{dx}\right\|_{H^{-1}} \leq C_2 \|\vphi\|_{L^2}
	\end{align*}
	for some $C_2 \geq 0$. Since those inequalities hold for all $\vphi \in C_c^\infty$, \eqref{e:hl_charac} follows for~$\ell = 1$ with appropriately defined $C \geq 0$.
	
	Conversely, suppose that \eqref{e:hl_charac} holds for all $\vphi\in C_c^\infty$. By~\Cref{l:dual_nu_dens}, this gives 
	\begin{align*}
		\left\Vert \calK_2^\ast (D^\alpha \vphi) \right\Vert_{H^{-1}} = \left\Vert\frac{d \nu_{D^\alpha \vphi}}{dx} \right\Vert_{H^{-1}} \leq C\|\vphi\|_{L^2}
	\end{align*}
	for every $|\alpha| = 1$. Hence, $|\<\calK f, D^\alpha \vphi \>| = | \langle f, K_2^\ast (D^\alpha \vphi) \rangle | \le C_\alpha\|f\|_{H^1}\|\vphi\|_{L^2}$ holds for all $f\in H^1$ and all $\vphi\in C_c^\infty$. Therefore, for fixed $f\in H^1$, the functional $\vphi\mapsto \<\calK f, D^\alpha \vphi \>$ is bounded on $C_c^\infty$ and can be extended to a bounded functional $F_f$ on $L^2$. By Riesz' representation theorem, $F_f\vphi = \<Lf,\vphi\>$ for all $\vphi\in L^2$ with some $Lf\in L^2$. Equating those expressions gives $\<\calK f,D^\alpha \vphi\> = F_f\vphi = \<Lf,\vphi\>$ for all $f\in H^1$ and all $\vphi\in C_c^\infty$. This shows that $\calK f$ has the weak derivative $-Lf = D^\alpha (\calK f) \in L^2$ for all $f\in H^1$ and all $|\alpha| = 1$, hence $\calK f\in H^1$ and $H^1$ is invariant under $\calK$. This shows the equivalence from (i) in \Cref{thm:hs_frac_charac}.
	
	\medskip
	\noindent  \textbf{Step 2: $H^\ell$-invariance, $\ell\in \bN$}. The characterization from step 1 is easily generalized to Sobolev spaces of higher smoothness. Fix $\ell \in \bN$. By definition, we have $f \in H^\ell$ if and only if $D^\alpha f \in L^2$ for every $\alpha \in \bN^n$ with $|\alpha| \leq \ell$. For the sufficiency, assume the invariance of $H^\ell$, i.e., $\calK f \in H^\ell$ for every $f \in H^\ell$. Again, this gives the boundedness of $\calK: H^\ell \to H^\ell$ by~\Cref{prop:RKHS_inv}, and for all $\alpha \in \bN^n$ with $|\alpha| \leq \ell$ we have $D^\alpha (\calK f) \in L^2$. By the definition of the weak derivative \eqref{eq:sobolev_weak_deriv_def}, we have 
	\begin{align*}
		\left| \left\< \calK f, D^\alpha \vphi \right\> \right| = \left| \left\< D^\alpha (\calK f), \vphi \right\> \right|
	\end{align*}
	for all $f \in H^\ell$ and $\vphi \in C_c^\infty$. Replacing $H^1$ by $H^\ell$ in the arguments from step 1 from \eqref{eq:h1_bound} onward, we get \eqref{e:hl_charac} for~$\ell \in \bN$. 
	
	Conversely, in analogy to step 1, \eqref{e:hl_charac} gives the boundedness of the functional $\vphi \mapsto \< \calK f, D^\alpha \vphi \>$ for fixed $f \in H^\ell$ and every $|\alpha| \leq \ell$ on $C_c^\infty$ and its extension on $L^2$. By Riesz' representation theorem, this functional is of form $\<Lf, \vphi \>$ with $Lf \in L^2$. We obtain the weak derivative $D^\alpha (\calK f) = (-1)^{|\alpha|} Lf$ for all $|\alpha| \leq \ell$ and hence the invariance of $H^\ell$ under $\calK$. 
	
	Thus, the invariance of $H^\ell$ under $\calK$ and \eqref{e:hl_charac} are equivalent, concluding the proof of part (i) of \Cref{thm:hs_frac_charac}.
	
	In the next step, we use an interpolation argument to deduce a sufficient condition for fractional Sobolev space invariance. 
	
	\medskip
	\noindent  \textbf{Step 3: $H^s$-invariance, $s\in (0,\ell]$}.
	By~\cite[Theorem 14.2.3]{BrenScot2008}, for real-valued $s \in (0, \ell]$, the space $H^s$ coincides with the real interpolation space
	\begin{align*}
		H^s(\Omega) = [L^2(\Omega), H^\ell(\Omega)]_\theta, \quad \theta = s/\ell \in (0,1),
	\end{align*}
	with equivalent norms. Under condition \eqref{e:L1}, we have the invariance of $L^2(\Omega)$ under $\calK$. Let \eqref{e:hl_charac} hold for~$\ell \in \bN$. Then $H^\ell(\Omega)$ is also invariant under $\calK$ by step 2. \cite[Theorem 22.3]{Tart2007} now implies that the Koopman operator
	\begin{align*}
		\calK: [L^2(\Omega), H^\ell(\Omega)]_\theta \to [L^2(\Omega), H^\ell(\Omega)]_\theta
	\end{align*}
	mapping from this interpolation space is well-defined and bounded. Thus, $H^s(\Omega)$ in indeed invariant under $\calK$, and part (ii) of \Cref{thm:hs_frac_charac} is proven. 
	
	The final two steps of the proof deal with the special case that the transition kernel is given by a differentiable transition density, in which case the Radon-Nikodym derivative $\frac{d \nu_g}{dx}$ can be explicitly calculated. 
	
	\medskip
	\noindent \textbf{Step 4: Deriving an explicit Radon-Nikodym density.} 
	Assume that the transition kernel $\rho_x$ admits a density w.r.t.\ the Lebesgue measure for every $x \in \Omega$, i.e., there exists some $p_x: \Omega \to [0, \infty)$ such that $\rho_x(\Delta) = \int_\Delta p_x(y) \,dy$ for every $\Delta \in \calB_\Omega$. Consider the family of densities as mapping $p: \Omega \times \Omega \to [0, \infty), p(x,y) = p_x(y)$ and write $d \rho_x(y) = p(x,y) \,d y$. This gives $\rho_x(\Delta) = \bP(Y \in \Delta \mid X = x) = \int_{\Omega} \chi{}_{\Delta}(y)  \,p(x,y) \, dy$ and for~$g \in L^\infty(\Omega)$, Fubini's theorem yields
	\begin{align*}
		\nu_g(\Delta) = \int_\Omega \biggl(\int_\Omega \chi_\Delta (y) p(x,y) \,dy \biggr) \, g(x) \, dx
		= \int_\Omega \chi_\Delta (y) \biggl(\int_\Omega p(x,y) g(x) \,dx \biggr) \, dy,
	\end{align*}
	so the Radon-Nikodym density $z_g = \frac{d \nu_g}{dx}$ of $\nu_g$ takes the form
	\begin{align}\label{eq:radon_nik_nu_smooth_dens}
		z_g(y) = \int_\Omega p(x,y) g(x) \, dx.
	\end{align}
	
	\medskip
	\noindent \textbf{Step 5: Sufficient condition via density derivatives.}
	To apply part (ii) of \Cref{thm:hs_frac_charac}, we take $g = D^\alpha \varphi$ with $\alpha \in \bN^n$, $|\alpha| \leq \ell \in \bN$, and $\varphi \in C_c^\infty$ in \eqref{eq:radon_nik_nu_smooth_dens}. Integration by parts yields the expression
	\begin{align*}
		z_{D^\alpha \vphi}(y) = \int_\Omega p(x,y) (D^\alpha\varphi)(x) \, dx = (-1)^{|\alpha|} \int_\Omega D_1^\alpha \,p(x,y) \varphi(x) \, dx
	\end{align*}
	for the Radon-Nikodym density of $\nu_{D^\alpha \vphi}$, where $D_1^\alpha$ denotes the $\alpha$-th weak derivative with respect to the first argument. Now, to compute $\left\| \frac{d \nu_{D^\alpha \vphi}}{dx} \right\|_{H^{-\ell}} = \| z_{D^\alpha \vphi} \|_{H^{-\ell}}$, we take $f \in H^\ell$ and calculate
	\begin{align*}
		\left\vert\left\langle f, z_{D^\alpha \vphi} \right\rangle_{H^{\ell}, H^{-\ell}}\right\vert &= \left\vert \int_\Omega f(y) z_{D^\alpha \vphi}(y) \,dy \right\vert \\
		&= \left\vert \int_\Omega f(y) \int_\Omega D_1^\alpha\,p(x,y) \vphi(x) \,dx \,dy \right\vert \\
		&\leq \Vert f\Vert_{L^2} \left\Vert \int_\Omega D_1^\alpha\, p(x,\cdot)\, \vphi(x)\,dx  \right\Vert_{L^2} \\&\leq \Vert f\Vert_{L^2} \Vert \vphi \Vert_{L^2} \Vert D_1^\alpha \,p(\cdot, \cdot) \Vert_{L^2(\Omega \times \Omega)} \leq \Vert f\Vert_{H^\ell}\, \Vert \vphi \Vert_{L^2} \Vert D_1^\alpha \,p(\cdot, \cdot) \Vert_{L^2(\Omega \times \Omega)}.
	\end{align*}
	In this, we applied the Cauchy-Schwarz inequality twice and then used that $\|f\|_{L^2} \leq \| f\|_{H^\ell}$. Recall that
	\begin{align*}
		\left\| z_{D^\alpha \vphi} \right\|_{H^{-\ell}} = \sup\left\{|\<f, z_{D^\alpha \vphi} \>_{H^\ell,H^{-\ell}}| : f\in H^\ell,\,\|f\|_{H^\ell}=1\right\}\, \le\, \Vert D_1^\alpha \,p(\cdot, \cdot) \Vert_{L^2(\Omega \times \Omega)} \,\|\vphi\|_{L^2}.
	\end{align*}
	Thus, if \eqref{e:Hs_condition_density} holds for some $C > 0$, we have
	\begin{align*}
		\sum_{|\alpha| \leq \ell} \left\| \frac{d \nu_{D^\alpha \vphi}}{dx} \right\|_{H^{-\ell}} \leq \sum_{|\alpha| \leq \ell} \Vert D_1^\alpha \,p(\cdot, \cdot) \Vert_{L^2(\Omega \times \Omega)} \,\|\vphi\|_{L^2} \leq C \,\|\vphi\|_{L^2}.
	\end{align*}
	By part (ii) of \Cref{thm:hs_frac_charac}, this implies that $H^s$ is invariant under $\calK$ for all $s \in (0, \ell]$. This is precisely the statement of part (iii) of \Cref{thm:hs_frac_charac}, concluding the proof. 
\end{proof}

\section{Error bounds for approximations of Koopman operators of stochastic systems}\label{ch:err_bounds}

This section is concerned with deriving error bounds for the approximation of the stochastic Koopman operator defined by
\begin{align}\label{eq:defkoop}
	(\calK f)(x) := \bE\big[f(Y)\,|\,X=x].
\end{align}
As our proof strategy reveals, we will rely on the following standing assumption.
\begin{assumption}\label{ass:rkhs_invariance}
	Let $k: \bR^n \times \bR^n \to \bR$ be a continuous strictly positive-definite symmetric kernel and $\bH$ the induced reproducing kernel Hilbert space on $\Omega$ with $\calK \in L(\bH,\bH)$ and ${\calK \in L(C_b(\Omega),C_b(\Omega))}$.
\end{assumption}

To verify the condition $\calK \in L(\bH,\bH)$ in~\Cref{ass:rkhs_invariance}, we refer to the sufficient conditions of the previous \Cref{ch:hs_invar}, where \Cref{thm:hs_frac_charac} guarantees the invariance of $\bH$ since Sobolev spaces $\bH = H^s$ are induced by Wendland kernels as RKHS \cite[Theorem 4.1]{KohnPhil24}. 

We start by constructing the kernel-based approximation $\wh \calK^\mathrm{MC}$ of the stochastic Koopman operator from \eqref{eq:koop_edmd_projection_stochastic} more thoroughly. In an intermediate step, we derive a kernel EDMD approximant $\wh \calK$, 
significantly building upon the analogue from \cite{KohnPhil24} in the deterministic setting. This allows to split the approximation error as follows:
\begin{align*}
	\|\calK - \wh\calK^{\mathrm{MC}}\|_{\bH \to C_b}
	\leq \underbrace{\|\calK - \wh\calK\|_{\bH \to C_b}}_{\text{\Cref{thm:err1} (det.)}} + \underbrace{\|\wh\calK - \wh\calK^{\mathrm{MC}}\|_{\bH \to C_b}}_{\text{\Cref{thm:mc_error1} (prob.)}}.
\end{align*}
The remainder of this section works towards the proof of \Cref{thm:total_error_bound} by analyzing both errors independently. The intermediary approximation $\wh \calK$ is employed in \Cref{subch:det_approx} to prove a bound on the deterministic projection error, see~\Cref{thm:err1}. Therein, we appeal to similar arguments as \cite{KohnPhil24}, where the invariance of a reproducing kernel Hilbert space under the Koopman operator is assumed as in our \Cref{ass:rkhs_invariance}. Since $\wh \calK$ requires the exact knowledge of propagated observables under the stochastic Koopman operator, i.e., conditional expectations, \Cref{subch:prob_approx} proposes the random, but computable approximation $\wh \calK^\mathrm{MC}$. The approximation of the ``exact'', but not computable $\wh \calK$ introduces a probabilistic error that is dealt with in \Cref{thm:mc_error1}.

\subsection{Deterministic approximation step: interpolation}\label{subch:det_approx}
The first (deterministic) approximation mimicks the approach of \cite{KohnPhil24}. However, this approximation will not be useful in practice as it involves expectation values and thus is not computable, in general. Therefore, we merely use it as a theoretical tool for an error bound of a computable estimator which we shall propose in the next subsection.

Let $\calX = \{x_1,\ldots,x_d\}\subset \Omega$ be a set of points. As described in \Cref{subch:stoch_kedmd}, we define a kernel-based approximation by
\begin{align}\label{eq:approx1}
	\widehat \calK = S_\calX \calK S_\calX,
\end{align}
where $S_\calX$ is the $\mathbb{H}$-orthogonal projection onto $V_\calX = \operatorname{span}\{\Phi_{x_i}: i\in\{1,\dots,d\}\} \subset \bH$ with the canonical features $\Phi_{x_i} = \kernel(x_i,\cdot)\in \mathbb{H}$. 
With the notation from \Cref{subch:stoch_kedmd}, recall the explicit projection formula \eqref{eq:projection_formula} given by $S_\calX g = g_\calX^\top \mathbf{K}_{\calX}^{-1} \mathbf{k}_\calX$ for some $g \in \bH$. The action of the abstract operator defined in \eqref{eq:approx1} now may explicitly be computed via
\[
\wh\calK f = S_\calX\calK S_\calX f = (\calK S_\calX f)_\calX^\top \mathbf{K}_{\calX}^{-1} \mathbf{k}_\calX = \sum_{i=1}^d \big[\mathbf{K}_{\calX}^{-1}\big(\underbrace{\mathbb{E}[(S_\calX f)(Y) \mid X=\cdot]}_{\calK S_\calX f}\big)_\calX\big]_i\Phi_{x_i}.
\]
Using linearity of the expected value and the projection formula \eqref{eq:projection_formula}, we get
\begin{align*}
	\mathbb{E}[S_\calX f (Y) \mid X=x_j] = f_\calX^\top \mathbf K_\calX^{-1} \mathbb{E}[\mathbf{k}_{\calX}(Y) \mid X=x_j] = \mathbb{E}[\mathbf{k}_{\calX}^\top(Y) \mid X=x_j] \mathbf K_\calX^{-1} f_\calX
\end{align*}
for each $j\in \{1,\dots,d\}$. Hence, defining the coefficients
\begin{equation*}
	\alpha_{\widehat \calK f} = \mathbf{K}_{\calX}^{-1}\bbK_{\calX, \calY}^\top\mathbf{K}_{\calX}^{-1}f_\calX, \qquad \mathrm{where}\qquad (\bbK_{\calX, \calY})_{i,j} = \mathbb{E}[\kernel(Y,x_i)\mid X = x_j],
\end{equation*}
yields the approximant
\begin{align}\label{eq:hatK_repr}
	\widehat\calK f = \sum_{i=1}^d (\alpha_{\widehat{K} f})_i \Phi_{x_i} = \alpha_{\widehat{K} f}^\top \mathbf k_\calX = f_\calX^\top \mathbf{K}_{\calX}^{-1}\bbK_{\calX, \calY} \mathbf{K}_{\calX}^{-1}\mathbf{k}_\calX.
\end{align}
In what follows, we shall assume that the probability kernel has the continuity-like property
\begin{align}\label{e:Cb}
	\forall\, x_0\in\Omega\;\forall\,\Delta\in\calB_\Omega(x_0) : \quad \rho_x(\Delta)\to\rho_{x_0}(\Delta)\quad\text{ as $x\to x_0$},
\end{align}
where $\calB_\Omega(x) = \{\Delta\in\calB_\Omega : \rho_x(\partial\Delta)=0\}$ for~$x\in\Omega$, with $\partial \Delta$ denoting the boundary of the Borel-set~$\Delta$. Note that this holds in the deterministic case where $\rho_x = \delta_{F(x)}$ with continuous $F: \Omega \to \Omega$: Then, $\rho_{x_0}(\partial \Delta) = \delta_{F(x_0)}(\partial \Delta) = 0$ is equivalent to $F(x_0) \notin \partial \Delta$, so $F(x_0)$ lies either in the interior or the exterior of $\Delta$, and continuity of $F$ gives the convergence. By the Portmanteau theorem \cite[Theorem 13.16]{Klen2008}, the condition \eqref{e:Cb} is equivalent to the fact that $\calK f$ is continuous for each $f\in C_b(\Omega)$. Hence, if \eqref{e:Cb} is satisfied, the Koopman operator is a contraction from $C_b(\Omega)$ into itself: for~$f\in C_b(\Omega)$, we have
\begin{align*}
	\|\calK f \|_\infty = \sup_{x \in \Omega} \left\vert \int_\Omega f \, d\rho_x \right\vert \leq \sup_{x \in \Omega} \int_\Omega \| f\|_\infty \, d\rho_x = \| f\|_\infty. 
\end{align*}
Combining the sufficient conditions of \Cref{thm:hs_frac_charac} with \eqref{e:Cb}, \Cref{ass:rkhs_invariance} holds, under which the approximant $\widehat\calK$ from \eqref{eq:hatK_repr} satisfies the following error bound.
\begin{theorem}[Deterministic error bound]\label{thm:err1}
	Under \Cref{ass:rkhs_invariance}, i.e. $\calK \in L(\bH, \bH)$ and $\calK\in L(C_b(\Omega),C_b(\Omega))$, the approximant $\wh \calK$ from \eqref{eq:hatK_repr} admits the error bound
	\begin{align}\label{eq:det_error}
		\big\|\calK - \widehat\calK\big\|_{\mathbb{H}\to C_b}\,\le\,
		(1+\|\calK\|_{\bH\to \bH})\cdot \|\calI - S_\calX\|_{\bH\to C_b}.
	\end{align}
\end{theorem}
\begin{proof}
	The proof is a particular case of the proof of \cite[Theorem 3.4]{KohnPhil24} and we repeat it here for the sake of clarity.
	As $\calK - \wh\calK = \calK - S_\calX\calK S_\calX = \calK (\calI - S_\calX) + (\calI - S_\calX)\calK S_\calX$, we get
	\begin{align*}
		\big\|\calK - \wh\calK\big\|_{\bH\to C_b}
		&\le \|\calK (\calI - S_\calX)\|_{\bH\to C_b} + \|(\calI - S_\calX)\calK S_\calX
		\|_{\bH\to C_b}\\
		&\le \|\calI - S_\calX\|_{\bH\to C_b}\,\left(\|\calK\|_{C_b\to C_b} + \|\calK\|_{\bH\to\bH}\|S_\calX\|_{\bH\to\bH}\right).
	\end{align*}
	Under the assumption $\calK \in L(C_b, C_b)$, $\calK$ is a contraction. Any constant $f \equiv c \in \Omega$ satisfies
	\begin{align*}
		\calK f(x) = \bE[f(Y) \mid X = x] = \bE [c \mid X = x] = c = f(x) \quad \forall x \in \Omega,
	\end{align*}
	so $\|\calK\|_{C_b\to C_b} = 1$. Since $S_\calX$ is a projection, we also have $\|S_\calX\|_{\bH\to\bH} = 1$, and the claim follows. 
\end{proof}

By~\Cref{prop:RKHS_inv}, the invariance assumption $\calK \in L(\bH, \bH)$ implies that the operator norm $\|\calK\|_{\bH\to\bH}$ is bounded. Thus, it is in order to focus on the projection error $\calI - S_\calX$ in dependence of the point set $\calX$. The projection error in \eqref{eq:det_error} attains a particularly intuitive form when considering Sobolev spaces as RKHS. As stated before, Wendland and Matérn kernels both generate such Sobolev spaces, allowing for the application of results of approximation theory (\cite[Theorem 11.13]{Wend04}) to further specify \Cref{thm:err1} in these cases. Under the appropriate sufficient condition from \Cref{ch:hs_invar}, the order of convergence can be quantified in terms of the smoothness of the underlying kernels. We exemplarily state the result for Wendland kernels below. 

\begin{corollary}[Deterministic error bound in Wendland native spaces]\label{cor:det_err_wendland}
	Let $l \in \bN$ and define $\sigma_{n,l} = \frac{n+1}{2}+l$. Consider the Wendland kernel $\kernel_{n,l}: \bR^n \times \bR^n \to \bR$ that generates the Sobolev space $H^{\sigma_{n,l}}(\Omega)$ as reproducing kernel Hilbert space $\bH$ with equivalence between the native norm $\|\cdot\|_\bH$ and the Sobolev norm $\|\cdot \|_{H^{\sigma_{n,l}}}$. Assume that $\calK \in L(C_b, C_b)$ and that \eqref{e:hl_charac} is satisfied for some $\ell \in \bN$ with $\ell > \sigma_{n,l}$. Then, the Koopman approximation $\wh \calK$ from \eqref{eq:approx1} satisfies
	\begin{align}\label{eq:det_err_wendland}
		\|\calK - \wh \calK \|_{\bH \to C_b} \leq C \,(1 +\|\calK\|_{\bH\to \bH})\cdot h_\calX^{l + 1/2}
	\end{align}
	with some constant $C > 0$ and the \emph{fill distance}
	\begin{align*}
		h_\calX = \sup_{x\in\Omega}\,\min_{i\in \{1,\dots,d\}}\|x - x_i\|_2.
	\end{align*}
\end{corollary}
\begin{proof}
	Under \eqref{e:hl_charac}, the Sobolev space $H^{\sigma_{n,l}}$ is invariant under the Koopman operator $\calK$. By~\cite[Theorem 11.17]{Wend04} or the restatement \cite[Theorem 5.1]{KohnPhil24}, we may bound the projection error by
	\begin{align*}
		\| \calI - S_\calX \|_{\bH \to C_b} \leq C h_\calX^{l + 1/2}. 
	\end{align*}
	\Cref{thm:err1} now yields \eqref{eq:det_err_wendland}. 
\end{proof}

\begin{remark}
	A similar result for Matérn kernels can be readily obtained; we state it in \Cref{app:matern_kernels}. 
	For equal smoothness parameters, the bound on the projection error in the case of Wendland kernels gains an extra $1/2$ in the exponent of the fill distance $h_\calX$ compared to Matérn kernels. This is because the locally supported Wendland kernels allow for a bound on the derivatives that scales with the fill distance, while the Matérn kernels are globally supported and stationary. In other words, the compact support introduces boundary effects that do not play a role for the global Matérn kernels. 
\end{remark}

For a large data set $\calX$, the kernel matrix $K_\calX$ is typically badly conditioned, leading to numerical instabilities. To alleviate this, one might regularize the kernel-based interpolation problem as in \cite{BoldPhil24} to get $R_\calX^\lambda f := {\arg \min}_{{g \in V_\calX}} \sum_{i=1}^d \left|f(x_i) - g(x_i)\right|^2 + \lambda \|g \|_\bH^2$
with regularization parameter $\lambda \geq 0$ and solution operator $R_\calX^\lambda$. The linear solution operator satisfies ${R_\calX^\lambda f = f_\calX^\top (\mathbf{K}_{\calX} + \lambda I)^{-1}\mathbf{k}_\calX}$ for~$f \in \bH$, so for~$\lambda = 0$ we recover $S_\calX$, while for~$\lambda > 0$, $R_\calX^\lambda: \bH \to \bH$ is no longer a projection. However, as is shown in \cite[Appendix A]{BoldPhil24}, in analogy to \eqref{eq:approx1} we can still define the approximant 
\begin{align}\label{eq:approx1_tikhonov}
	\widehat{\calK}_\lambda := R_\calX^\lambda \,\calK \,R_\calX^\lambda
\end{align}
that acts on $f \in \bH$ via
\begin{align}\label{eq:hatK_tichonov}
	\widehat{\calK}_\lambda f = f_\calX^\top \,(\mathbf{K}_{\calX} + \lambda I)^{-1}\,\bbK_{\calX, \calY}\, (\mathbf{K}_{\calX} + \lambda I)^{-1}\,\mathbf{k}_\calX.  
\end{align}
For this \emph{regularized} approximant, we get a similar error bound as in \Cref{thm:err1}.
\begin{proposition}[Regularized deterministic error bound]\label{prop:err_tikh}
	Under \Cref{ass:rkhs_invariance}, the regularized approximant \eqref{eq:hatK_tichonov} with $\lambda \geq 0$ satisfies
	\begin{align}\label{eq:det_error_tikhonov}
		\bigl\Vert \calK - \widehat{\calK}_\lambda \bigr\Vert_{\bH\to C_b} \leq \left( 1 + \| \calK \|_{\bH \to \bH}  \right)\, \| \calI - R_\calX^\lambda\|_{\bH\to C_b}.
	\end{align}
	Specifically, for the RKHS $\bH = H^{\sigma_{n,l}}(\Omega)$ generated by the Wendland kernel $\kernel_{n,l}$, where $\sigma_{n,l} = \frac{n+1}{2}+l$, the deterministic error admits the bound
	\begin{align}\label{eq:det_error_tikh_wendland}
		\| (\calK - \wh \calK_\lambda) f \|_\infty \leq C \left( h_\calX^{l+1/2} + \lambda^{1/2} \right) \| f \|_\bH.
	\end{align}
\end{proposition}
\begin{proof}
	The proof of \eqref{eq:det_error_tikhonov} follows the lines of the proof of \Cref{thm:err1} and uses the fact that $\bigl\Vert R_\calX^\lambda \bigr\Vert_{\bH \to \bH} \leq 1$ as shown in \cite[Proposition B.1 (iii)]{BoldPhil24}. The special case covered in \eqref{eq:det_error_tikh_wendland} follows from \cite[Theorem 2.4]{BoldPhil24}.
\end{proof}
The second part of \Cref{prop:err_tikh} indicates the bias introduced by the regularization parameter $\lambda$. For decreasing fill distance, the regularization bias becomes dominant in the error bound. As a tradeoff, regularization allows for a control on the probabilistic part of the error by bounding the smallest eigenvalue of $\bbK_\calX$ from below, see~\Cref{cor:mc_error1_tikh} in the following section. 

\subsection{Probabilistic approximation step: estimating the expected value}\label{subch:prob_approx}

In practice, we cannot form $\bbK_{\calX, \calY}$, as we do not have access to the expected values. For this reason, we introduce a random matrix ${\mathbf{K}}^{\mathrm{MC}}_{\calX, \calY}$, whose entries are Monte Carlo approximations of the expected values in~$\bbK_{\calX, \calY}$. We denote the $j$-th column of $\bbK_{\calX, \calY}$ by $(\bbK_{\calX, \calY})_j$. The $j$-th column of $\bbK_{\calX, \calY}^\mathrm{MC}$ is then defined by
\begin{align}\label{eq:MC}
	\big(\mathbf{K}^{\mathrm{MC}}_{\calX, \calY}\big)_j:= \frac{1}{m_j} \sum_{l=1}^{m_j} \bbk_\calX(y_{l,j}) \approx \mathbb{E}[\bbk_\calX(Y)\mid X=x_j],
\end{align}
where $y_{l,j}$, $l\in \{1,\dots,m_j\}$, are drawn i.i.d.\ from $\rho(x_j,\cdot)$, $j\in \{1,\dots,d\}$.

This leads to a \emph{computable} approximation $\widehat \calK^{\mathrm{MC}}$, which is defined analogously to $\widehat\calK$ in \eqref{eq:hatK_repr}, with the only difference that we replace $\bbK_{\calX, \calY}$ with $\mathbf{K}^\mathrm{MC}_{\calX, \calY}$ as defined in \eqref{eq:MC}. Thus, we set
\begin{align}\label{eq:mcedmd}
	\widehat \calK^{\mathrm{MC}}f = f_\calX^\top \,\mathbf{K}_{\calX}^{-1}\,\bbK_{\calX, \calY}^{\mathrm{MC}}\, \mathbf{K}_{\calX}^{-1}\,\bbk_\calX.
\end{align}
Decomposing the \emph{full error} by inserting the auxiliary approximation $\widehat \calK$, i.e.,  
\begin{align*}
	\|\calK - \wh\calK^{\mathrm{MC}}\|_{\bH \to C_b}
	\leq \underbrace{\|\calK - \wh\calK\|_{\bH \to C_b}}_{\text{\Cref{thm:err1} (det.)}} + \underbrace{\|\wh\calK - \wh\calK^{\mathrm{MC}}\|_{\bH \to C_b}}_{\text{Monte Carlo (prob.)}},
\end{align*}
we are left with bounding the second term, after analyzing the deterministic first term in \Cref{subch:det_approx}. Due to the random sampling employed in~$\bbK_{\calX, \calY}^\mathrm{MC}$, this error is probabilistic. In the remainder of this part, we derive a probabilistic uniform bound on the error
\begin{align*}
	(\wh\calK - \wh\calK^{\mathrm{MC}})f = f_\calX^\top\, \mathbf{K}_{\calX}^{-1}\,\left(\bbK_{\calX, \calY} - \bbK_{\calX, \calY}^{\mathrm{MC}}\right)\, \mathbf{K}_{\calX}^{-1} \,\mathbf{k}_\calX.
\end{align*}
For this, let us define the vectors
\begin{align*}
	g^\top := f_\calX^\top \mathbf{K}_{\calX}^{-1}\bbK_{\calX, \calY},
	\qquad
	(g^\mathrm{MC})^\top= f_\calX^\top \mathbf{K}_{\calX}^{-1}\bbK_{\calX, \calY}^{\mathrm{MC}}.
\end{align*}
Via Hölder's inequality, we rewrite the probabilistic error term in terms of those vectors. Subsequently, we bound the maximum norm of the Monte Carlo estimation error between $g$ and $g^\mathrm{MC}$ via the Hoeffding inequality. 
\begin{theorem}[Monte Carlo error bound]\label{thm:mc_error1}
	Choose $m_j = m \in \bN$ i.i.d.\ samples from every $\rho_{x_j}$, $j = 1, \dots, d$, and let $\varepsilon > 0$ such that $m \geq \frac{2\|\kernel\|_\infty}{\veps^2}\log(2d)$. Then, for all $f \in \bH$, the Koopman approximation $\wh\calK^{\mathrm{MC}}$ from \eqref{eq:mcedmd} satisfies
	\begin{align}\label{eq:MC_error}
		\left\|(\wh\calK^{\mathrm{MC}} - \wh\calK)f\right\|_\infty\,\le\,\veps\|\kernel\|_\infty^{1/2}\cdot\sqrt{\max_{v\in\one}v^\top \mathbf{K}_\calX^{-1}v}\cdot\|f\|_\bH
	\end{align}
	with probability at least $1 - 2d \exp\left(-\frac{m\veps^2}{2\|\kernel\|_\infty}\right)$, where $\one = \{v\in\R^d : |v_i|=1\}$.
\end{theorem}
\begin{proof}
	For all $x\in \Omega$, we may interpret this as a perturbed interpolation problem, where the perturbation of the function values stems from approximation of the expected value: 
	\begin{align*}
		\left|(\widehat \calK^{\mathrm{MC}}-\widehat \calK)f(x)\right| = \left|(g - g^\mathrm{MC})^\top \bbK_\calX^{-1} \mathbf{k}_\calX(x)\right| \leq \left\|g-g^{\mathrm{MC}}\right\|_\infty \, \left\|\bbK_\calX^{-1} \mathbf{k}_\calX(x)\right\|_1.
	\end{align*}
	As shown in the proof of \cite[Theorem 4.2]{BoldPhil24}, the deterministic second factor is bounded by $\|\bbK_\calX^{-1} \mathbf{k}_\calX(x)\|_1 \le \|\kernel\|_\infty^{1/2}\cdot\big(\max_{v\in\one}v^\top \mathbf{K}_\calX^{-1}v\big)^{1/2}$.
	Moreover,
	\begin{align*}
		\|g-g^\mathrm{MC}\|_\infty = 
		\max_{x_j\in \calX}\left|\mathbb{E}[(S_\calX f)(Y)\mid X=x_j] - \frac{1}{m_j} \sum_{l=1}^{m_j} (S_\calX f)(y_{l,j})\right|.
	\end{align*}
	Setting $\bE_j := g_j = \mathbb{E}[(S_\calX f)(Y)\mid X=x_j]$, we obtain $\bE[(S_\calX f)(y_{l,j})] = \int (S_\calX f)(y)\,\rho(x_j,dy) = (\calK S_\calX f)(x_j) = \bE_j$. We also have
	\[
	|(S_\calX f)(y)| = |\<S_\calX f,\Phi_y\>_\bH|\,\le\,\|S_\calX f\|_H\|\Phi_y\|_\bH\,\le\,\|\kernel\|_\infty^{1/2}\|f\|_\bH
	\]
	for all $y\in\Omega$, where the last inequality follows from $\| S_\calX \|_{H \to H} = 1$ since $S_\calX$ is an orthogonal projection. 
	Hence, the random variables $(S_\calX f)(y_{l,j})$ are bounded and i.i.d., so that we may apply Hoeffding's inequality:
	\begin{align*}
		\bP\left(\left|\bE_j - \frac{1}{m_j} \sum_{l=1}^{m_j} (S_\calX f)(y_{l,j})\right|\,\ge\,\veps\|f\|_H\right)
		&= \bP\left(\left|\sum_{l=1}^{m_j} \big[(S_\calX f)(y_{l,j}) - \bE_j\big]\right|\,\ge\,\veps m_j\|f\|_H\right)\\
		&\le 2\exp\left(-\frac{2m_j^2\veps^2\|f\|_H^2}{4m_j\|\kernel\|_\infty\|f\|_H^2}\right) = 2\exp\left(-\frac{m_j\veps^2}{2\|\kernel\|_\infty}\right).
	\end{align*}
	Plugging this entry-wise result into a union bound, we readily obtain:
	\begin{align*}
		\bP\big(\|g-g^\mathrm{MC}\|_\infty\ge\veps\|f\|_H)
		&= \bP\left(\exists j\in \{1,\dots,d\} : \left|\bE_j - \frac{1}{m_j} \sum_{l=1}^{m_j} (S_\calX f)(y_{l,j})\right|\,\ge\,\veps\|f\|_H\right)\\
		&\le \sum_{j=1}^d \bP\left(\left|\bE_j - \frac{1}{m_j} \sum_{l=1}^{m_j} (S_\calX f)(y_{l,j})\right|\,\ge\,\veps\|f\|_H\right)\\
		&\le 2\sum_{j=1}^d \exp\left(-\frac{m_j\veps^2}{2\|\kernel\|_\infty}\right).
	\end{align*}
	Hence, if $m_j = m$ for all $j\in \{1,\dots,d\}$, then
	\begin{align*}
		\bP\big(\|g-g^\mathrm{MC}\|_\infty\ge\veps\|f\|_H\big)
		&\le 2d\exp\left(-\frac{m\veps^2}{2\|\kernel\|_\infty}\right).
	\end{align*}
	In total, for every $f\in H$, we obtain the upper bound
	\begin{align*}
		\|(\wh\calK^{\mathrm{MC}} - \wh\calK)f\|_\infty\,\le\,\veps\cdot\|\kernel\|_\infty^{1/2}\cdot\Big(\max_{v\in\one}v^\top \mathbf{K}_\calX^{-1}v\Big)^{1/2}\|f\|_H,
	\end{align*}
	i.e., \eqref{eq:MC_error}, with probability at least $1 - 2d \exp\left(-\frac{m\veps^2}{2\|\kernel\|_\infty}\right)$, which is nonnegative due to the imposed inequality $m \geq \frac{2\|\kernel\|_\infty}{\veps^2}\log(2d)$.
\end{proof}
Note that in the bound \eqref{eq:MC_error} we have $\|\kernel\|_\infty < \infty$ by~\Cref{ass:rkhs_invariance}, since the kernel $\kernel$ is continuous on $\bR^n \times \bR^n$ and hence bounded on the bounded open set $\Omega \times \Omega$.  
\begin{remark}[Alternative probabilistic bound] \label{rem:l2_l2_mc_error_bound}
	Alternatively, one might use Hölder's inequality with $p = q = 2$:
	\begin{align*}
		|(\widehat \calK^{\mathrm{MC}}-\widehat \calK)f(x)| \leq \|g-g^{\mathrm{MC}}\|_2 \, \|K_\calX^{-1} \mathbf{k}_\calX(x)\|_2.
	\end{align*}
	In contrast to \Cref{thm:mc_error1}, which required bounding the maximum norm, this approach necessitates a bound on $\Vert g - g^\mathrm{MC}\Vert_2$ with high probability. Applying the Hilbert-space Hoeffding inequality circumvents the union bound used in the proof of \Cref{thm:mc_error1} to directly get 
	\begin{align}\label{eq:l2_l2_err_prob}
		\bP\big(\|g-g^\mathrm{MC}\|_2\ge\veps\|f\|_H) \leq 2\exp\left(-\frac{m\veps^2}{8d\|\kernel\|_\infty}\right).
	\end{align}
	By the $\ell^1$-$\ell^2$-inequality, the deterministic second term is smaller than the one in \Cref{thm:mc_error1}, i.e., $\|\bbK_\calX^{-1}\bbk_\calX(x)\|_2 \leq \|\bbK_\calX^{-1}\bbk_\calX(x)\|_1$. We refer the reader to \Cref{app:l2_l2_MC} for detailed calculations. In total, with probability at least $1 - 2\exp\left(-\frac{m\veps^2}{8d\|\kernel\|_\infty}\right)$, the bound
	\begin{align*}
		\|(\wh\calK^{\mathrm{MC}} - \wh\calK)f\|_\infty\,\le\,\veps\cdot\|\bbK_\calX^{-1}\bbk_\calX(x)\|_2\,\|f\|_H.
	\end{align*}
	holds. However, comparing the error probabilities, we notice that in \Cref{thm:mc_error1} the factor $d$ appears in front of the exponential term, while in \eqref{eq:l2_l2_err_prob}, the factor $d$ appears \emph{inside} the exponential term. From a practical standpoint, the former is preferable, since we only have to increase the number of Monte Carlo samples $m$ by a logarithmic factor in~$d$ to achieve the same error probability. In the latter case, the growth of $m$ would have to be linear in~$d$. 
\end{remark}

\begin{remark}[Impact of kernel matrix eigenvalues]\label{rem:kernel_matrix_cond}
	The bound
	\begin{align*}
		\|\bbK_\calX^{-1} \mathbf{k}_\calX(x)\|_1 \le \|\kernel\|_\infty^{1/2}\cdot\sqrt{\max_{v\in\one}v^\top \mathbf{K}_\calX^{-1}v}
	\end{align*}
	from the proof of \Cref{thm:mc_error1} warrants further discussion. While the second term might prove sharper, a more explicit upper bound involving the Rayleigh quotient can be calculated in the following way: As a first step, apply the $\ell^1$-$\ell^2$-norm inequality to get
	\begin{align*}
		\|\bbK_\calX^{-1} \mathbf{k}_\calX(x)\|_1^2 \leq d \cdot \|\bbK_\calX^{-1} \mathbf{k}_\calX(x)\|_2^2 &= d \cdot \bbk_\calX(x)^\top \bbK_\calX^{-2}\, \bbk_\calX(x) \\
		&\leq d \cdot \lambda_\mathrm{max}(\bbK_\calX^{-1})\,\bbk_\calX(x)^\top \bbK_\calX^{-1} \,\bbk_\calX(x).
	\end{align*}
	From the positive definiteness of the kernel $\kernel$, we get the positive semidefiniteness of the augmented kernel matrix $
	\bbK_{\calX \cup \{x\}} = \left(\begin{smallmatrix}
		\bbK_\calX & \bbk_\calX \\
		\bbk_\calX^\top & \kernel(x,x)
	\end{smallmatrix}\right)$
	with positive semidefinite (thus, non-negative) Schur complement $\kernel(x,x) - \bbk_\calX^\top \bbK_\calX^{-1} \bbk_\calX \geq 0$. With $\lambda_\mathrm{max}(\bbK_\calX^{-1})= \lambda_\mathrm{min}(\bbK_\calX)^{-1}$,  where $\lambda_{\max}$ and $\lambda_{\min}$ are the largest and smallest eigenvalue of the respective matrices, we get
	\begin{align*}
		\|\bbK_\calX^{-1} \mathbf{k}_\calX(x)\|_1^2 \leq \kernel(x,x) \cdot \frac{d}{\lambda_\mathrm{min}(\bbK_\calX)}.
	\end{align*}
	The term in the denominator usually decays rapidly as the condition of the kernel matrix gets worse. If our kernel is given by a \emph{radial basis function} (RBF) $B: \bR^n \to \bR$, $\kernel(x,x) = B(0)$ is constant, so 
	\begin{align*}
		\| \mathbf{K}_\calX^{-1} \mathbf{k}_\calX(x) \|_1^2 \leq  B(0) \frac{d}{\lambda_\mathrm{min}(\mathbf{K}_\calX)}.
	\end{align*} 
	For example, for the Wendland kernel $\kernel_{n,l}$ with RBF $B_{n,l}$, we have $B_{n,l}(0) = \frac{1}{2^l} \binom{\lfloor n/2\rfloor + 2l + 1}{l}$, while for Matérn kernels $\kernel_\nu$ with RBF $\calM_{\nu}$, we have $\calM_{\nu}(0) = 1$. 
\end{remark}

To balance the decay of the eigenvalue $\lambda_\mathrm{min}(\mathbf{K}_\calX)$ to zero, the computable Monte Carlo approximation \eqref{eq:mcedmd} also lends itself to a Tikhonov regularization in style of \eqref{eq:approx1_tikhonov}. Hence, we define the computable regularized approximation
\begin{align}\label{eq:regularized_mcemdmd}
	\widehat{\calK}_\lambda^\mathrm{MC} f = f_\calX^\top\, (K_\calX + \lambda I)^{-1}\, \bbK_{\calX, \calY}^\mathrm{MC}\, (K_\calX + \lambda I)^{-1} \,\mathbf{k}_\calX.
\end{align}
A similar error bound can be computed for the Monte Carlo estimation error $\widehat{\calK}_\lambda^\mathrm{MC} - \widehat{\calK}_\lambda$ in the regularized case. For this, we define 
\begin{align*}
	g_\lambda^\top := f_\calX^\top (\mathbf{K}_{\calX} + \lambda I)^{-1}\bbK_{\calX, \calY},
	\qquad
	(g_\lambda^\mathrm{MC})^\top= f_\calX^\top (\mathbf{K}_{\calX} + \lambda I)^{-1}\bbK_{\calX, \calY}^{\mathrm{MC}}.
\end{align*}
In similar fashion to \Cref{thm:mc_error1}, we obtain the following error bound. 
\begin{corollary}[Regularized Monte Carlo error bound]\label{cor:mc_error1_tikh}
	Let the assumptions of \Cref{thm:mc_error1} hold. For~$f \in \bH$, the regularized approximant $\wh\calK_\lambda^{\mathrm{MC}}$ from \eqref{eq:regularized_mcemdmd} satisfies
	\begin{align}\label{eq:mc_error_regularized}
		\left\|(\wh\calK_\lambda^{\mathrm{MC}} - \wh\calK_\lambda)f\right\|_\infty\,\le\,\veps\cdot\|\kernel\|_\infty^{1/2}\cdot \sqrt{\frac{d}{\lambda}}\, \|f\|_{\bH}
	\end{align}
	with probability at least $1 - 2d \exp\left(-\frac{m\varepsilon^2}{2 \|\kernel\|_\infty}\right)$.
\end{corollary}
\begin{proof}
	Again, we factor out the probabilistic part of the Monte Carlo error:
	\begin{align*}
		|(\widehat \calK_\lambda^{\mathrm{MC}}-\widehat \calK_\lambda)f(x)| &= |(g_\lambda - g_\lambda^\mathrm{MC})^\top (K_\calX + \lambda I)^{-1} \mathbf{k}_\calX(x)| \\
		&\leq \|g_\lambda-g_\lambda^{\mathrm{MC}}\|_\infty \| (K_\calX + \lambda I)^{-1}\mathbf{k}_\calX(x)\|_1.
	\end{align*}
	Regarding the deterministic second part of the error, from the positive definiteness of $\bbK_\calX$, we get $\lambda_{\min}(\bbK_\calX + \lambda I) = \lambda_{\min}(\bbK_\calX) + \lambda \geq \lambda$. Following \Cref{rem:kernel_matrix_cond}, this yields
	\begin{align*}
		\| (\mathbf{K}_\calX + \lambda I)^{-1} \mathbf{k}_\calX(x) \|_1^2 \leq \kernel(x,x) \frac{d}{\lambda_\mathrm{min}(\mathbf{K}_\calX + \lambda I)} \leq \kernel(x,x) \frac{d}{\lambda}.
	\end{align*}
	For the other term, we have
	\begin{align*}
		\|g_\lambda-g_\lambda^\mathrm{MC}\|_\infty = \max_{x_j\in \calX}\left|\mathbb{E}[(R_\calX^\lambda f)(Y)\mid X=x_j] - \frac{1}{m_j} \sum_{l=1}^{m_j} (R_\calX^\lambda f)(y_{l,j})\right|
	\end{align*}    
	with the operator $R_\calX^\lambda$ from \eqref{eq:hatK_tichonov}. By~\cite[Proposition B.1]{BoldPhil24}, that ``projection-like'' operator satisfies $\| R_\calX^\lambda \|_{\bH \to \bH} \leq \| S_\calX \|_{\bH \to \bH} = 1$. Thus, the arguments from the proof of \Cref{thm:mc_error1} carry over, and for every $f \in \bH$, we get
	\begin{align*}
		\bP\left(\|g_\lambda-g_\lambda^\mathrm{MC}\|_\infty \geq \veps \| f \|_\bH \right) \leq 2d\exp\left( -\frac{m\veps^2}{2\|\kernel\|_\infty} \right).
	\end{align*}
	In total, this yields \eqref{eq:mc_error_regularized}.
\end{proof}
The appeal of \Cref{cor:mc_error1_tikh} lies in the more explicit bound in terms of the number $d$ of data points in~$\calX$. The regularization parameter $\lambda$ puts a ``floor'' on the smallest eigenvalue of the kernel matrix $\bbK_\calX$, making the tradeoff between deterministic regularization bias and kernel matrix conditioning evident. 

In the final subsection, we combine the deterministic and Monte Carlo error bounds to get a bound on the full $L^\infty$-error of the approximations $\wh\calK^{\mathrm{MC}}$ and $\wh\calK_\lambda^{\mathrm{MC}}$, proving the remaining main results from \Cref{ch:setting}.

\subsection{\texorpdfstring{Proof of \Cref{thm:total_error_bound,thm:total_err_wendland}}{Proof of Theorem 2.6 and 2.8}}\label{subch:proof_total_approx}

In this part, we combine the results of \Cref{thm:err1,thm:mc_error1} to obtain the full error bound from \Cref{thm:total_error_bound}. After the proof, we provide the respective result for the regularized approximant $\widehat{\calK}_{\lambda}^{\mathrm{MC}}$ and a more explicit version in the central special case of a Sobolev space $\bH = H^s$ generated by a Wendland kernel. First, let us restate our main result for the stochastic kEDMD approximant.

\restatetheorem[$L^\infty$-error bound, restated]{thm:total_error_bound}
Let $\bH$ be an RKHS of functions over $\Omega$, and assume $\calK\in L(\mathbb{H}, \mathbb{H})$ and $\calK\in L(C_b(\Omega),C_b(\Omega))$. Furthermore, for~$\mathcal{X} = \{x_1,\ldots,x_d\}$, let $m_j = m \in \bN$ i.i.d.\ samples from every $\rho_{x_j}$, $j \in \{1, \dots, d\}$, be given.
Then, for every $\veps > 0$ satisfying $m \geq \frac{2\|\kernel\|_\infty}{\veps^2}\log(2d)$ and  for all $f \in \bH$, $\| f \|_\bH \leq 1$, we have the pointwise error bound
\begin{equation}
	\big\| \big(\calK - \widehat{\calK}^{\mathrm{MC}}\big) f\big\|_{\infty} \leq \big[ 1 + \|\calK\|_{\bH\to \bH} \big] \|\calI - S_\calX \|_{\bH \to C_b} + \veps \cdot \|\kernel\|_\infty^{1/2} \sqrt{\max_{v\in\one}v^\top \mathbf{K}_\calX^{-1}v} \nonumber
\end{equation}
with probability at least $1 - 2d \exp\left(-\frac{m\veps^2}{2\|\kernel\|_\infty}\right)$, where~$\one = \{v\in\R^d : |v_i|=1\}$.
\restatetheoremend

\begin{proof}
	Let $f \in \bH$ with $\|f\|_\bH \leq 1$. We decompose the full approximation error into its deterministic and probabilistic part by
	\begin{align*}
		\| (\calK - \wh\calK^{\mathrm{MC}})f\|_{\infty} \leq \| (\calK - \wh\calK)f\|_{\infty} + \| (\wh\calK - \wh\calK^{\mathrm{MC}})f\|_{\infty}.
	\end{align*}
	As $\calK \in L(\bH, \bH)$ and $\calK \in L(C_b, C_b)$, \Cref{ass:rkhs_invariance} is fulfilled. \Cref{thm:err1} now bounds the first summand by
	\begin{align*}
		\big\|(\calK - \widehat\calK\big) f\|_{\mathbb{H}\to C_b}\,\le\, (1+\|\calK\|_{\bH\to \bH})\cdot \|\calI - S_\calX\|_{\bH\to C_b} \| f \|_\bH.
	\end{align*}
	By~\Cref{thm:mc_error1}, the second summand admits the bound
	\begin{align*}
		\|(\wh\calK^{\mathrm{MC}} - \wh\calK)f\|_\infty\,\le\,\veps\cdot\|\kernel\|_\infty^{1/2}\cdot\sqrt{\max_{v\in\one}v^\top \mathbf{K}_\calX^{-1}v}\,\|f\|_H
	\end{align*}
	with probability at least $1 - 2d \exp\left(-\frac{m\veps^2}{2\|\kernel\|_\infty}\right)$. With $\|f\|_\bH \leq 1$, those two inequalities yield the claim.
\end{proof}

The projection error $\calI - S_\calX$ and the expression involving the kernel matrix $\mathbf{K}_\calX^{-1}$ warrant a closer inspection in certain special cases. For Sobolev spaces as Wendland RKHS, the projection error in \Cref{thm:total_err_wendland} is bounded in terms of the fill distance $h_\calX$.

\restatetheorem[$L^\infty$-error bound for Wendland native spaces, restated]{thm:total_err_wendland}
Let $l \in \bN$, define $\sigma_{n,l} = \frac{n+1}{2}+l$ and denote by $\bH$ the Sobolev space $H^{\sigma_{n,l}}(\Omega)$ induced by the Wendland kernel $\kernel_{n,l}$ as RKHS. Assume that $\calK \in L(\bH, \bH)$ and $\calK \in L(C_b(\Omega), C_b(\Omega))$. Furthermore, for~$\mathcal{X} = \{x_1,\ldots,x_d\} \subset \Omega$, let $m_j = m \in \bN$ i.i.d.\ samples from every $\rho_{x_j}$, $j \in \{1, \dots, d\}$, be given.
Then, for every $\veps > 0$ satisfying $m \geq \frac{2 C_2^2}{\veps^2}\log(2d)$ and all $f \in \bH$ with $\|f\|_\bH \leq 1$, we have the pointwise error bound
\begin{equation}
	\nonumber
	\begin{aligned}
		\| (\calK - \calK^{\mathrm{MC}})f\|_{\infty} \leq\, &C_1\,(1+\|\calK\|_{\bH\to \bH})\cdot h_\calX^{l + 1/2} + \veps \cdot C_2\cdot\sqrt{\max_{v\in\one}v^\top \mathbf{K}_\calX^{-1}v}
	\end{aligned}
\end{equation}
with probability at least $1 - 2d\exp\left(-\frac{m\veps^2}{2 C_2^2}\right)$, where $C_1, C_2 > 0$ are constants only depending on the kernel $\kernel_{n,l}$, and $\one = \{v\in\R^d : |v_i|=1\}$.
\restatetheoremend
\begin{proof}
	The result follows analogously to the proof of \Cref{thm:total_error_bound} by the error decomposition
	\begin{align*}
		\| (\calK - \wh\calK^{\mathrm{MC}})f\|_{\infty} \leq \| (\calK - \wh\calK)f\|_{\infty} + \| (\wh\calK - \wh\calK^{\mathrm{MC}})f\|_{\infty}.
	\end{align*}
	for~$f \in \bH$. By~\Cref{cor:det_err_wendland}, 
	\begin{align*}
		\|\calK - \wh \calK \|_{\bH \to C_b} \leq C_1 (1 +\|\calK\|_{\bH\to \bH})\cdot h_\calX^{l + 1/2}
	\end{align*}
	for some $C_1 > 0$. The second summand is bounded by~\Cref{thm:mc_error1}, which gives
	\begin{align*}
		\|(\wh\calK^{\mathrm{MC}} - \wh\calK)f\|_\infty\,\le\,\veps\cdot\|\kernel_{n,l}\|_\infty^{1/2}\cdot\sqrt{\max_{v\in\one}v^\top \mathbf{K}_\calX^{-1}v}\|f\|_H
	\end{align*}
	with probability at least $1 - 2d\exp\left(-\frac{m\veps^2}{2 \|\kernel_{n,l}\|_\infty}\right)$. Setting $C_2 = \|\kernel_{n,l}\|_\infty^{1/2}$ and taking both inequalities for~$\|f\|_\bH \leq 1$ yields the claim.
\end{proof}

Note that Wendland kernels are defined via $\kernel_{n,l}(x,y) = B_{n,l}(\|x-y\|_2)$ with the radial basis function $B_{n,l}: [0, \infty) \to \bR$, see~\Cref{app:wendland_kernels}. It can be shown that $\|\kernel_{n,l} \|_\infty = \sup_{x \in \bR^n} |\kernel_{n,l}(x,x)| = B_{n,l}(0) = p_{n,l}(0)$ (see, e.g., \cite{PhilSchaWort2025}), where $p_{n,l}$ is the recursively defined polynomial underlying the Wendland kernel \cite{Wend04}. Thus, we get the explicitly computable expression $C_2 = \sqrt{p_{n,l}(0)}$ for the constant from \eqref{eq:total_err_wendland} that depends on dimension $n$ and smoothness $l$. A similar bound as in \Cref{thm:total_err_wendland} is easily proven for Matérn kernels, which we state in \Cref{cor:total_error_bound_matern} in \Cref{app:matern_kernels}. 

To make use of the error bounds from \Cref{thm:total_err_wendland} in practice, one would first choose a large enough $d$ to ensure that the fill distance $h_\calX$ controlling the deterministic projection error is sufficiently small. This might result in a badly conditioned kernel matrix $\mathbf{K}_\calX$, making the expression $\sqrt{\max_{v\in\one}v^\top \mathbf{K}_\calX^{-1}v}$ large. 
Thus, to control the probabilistic error, one must choose a small enough $\veps > 0$ that compensates for the conditioning. Lastly, we want our error bound to hold with probability close to one, so the amount of Monte Carlo samples $m$ for the entries of $\bbK_{\calX, \calY}^\mathrm{MC}$ has to guarantee that $2d \exp(-\frac{m \veps^2}{2 \|\kernel\|_\infty})$ is close to zero. Since the factor involving the inverse kernel matrix $\bbK_\calX^{-1}$ depends on the size $d$ of data set $\calX$, regularization 
as indicated in \Cref{cor:mc_error1_tikh} prevents this expression from blowing up with $d$. However, regularization introduces a bias to the deterministic error, which is why the final result of this section explicitly juxtaposes deterministic and probabilistic error for the regularized approximant $\wh \calK_\lambda^\mathrm{MC}$. 
\begin{corollary}[Regularized $L^\infty$-error bounds]\label{cor:reg_total_err_bound}
	Let the assumptions of \Cref{thm:total_error_bound} hold and let $\lambda \geq 0$. Then, for every $\veps > 0$ and every $f \in \bH$ with $\|f\|_\bH \leq 1$, we have 
	\begin{equation}\label{eq:reg_total_err}
		\begin{aligned}
			\big\| \big(\calK - \wh \calK_\lambda^\mathrm{MC}\big) f \big\|_\infty \leq \left( \| 1 + \| \calK \|_{\bH \to \bH}  \right)\, \| \calI - R_\calX^\lambda\|_{\bH\to C_b} \,+ \veps\cdot\|\kernel\|_\infty^{1/2}\cdot \sqrt{\frac{d}{\lambda}}
		\end{aligned}
	\end{equation}
	with probability at least $1 - 2 d\exp\left( -\frac{m\varepsilon^2}{2 \|\kernel\|_\infty}\right)$. Specifically, for the RKHS $\bH = H^{\sigma_{n,l}}(\Omega)$ induced by the Wendland kernel $\kernel_{n,l}$, where $\sigma_{n,l} = \frac{n+1}{2}+l$, we have the pointwise error bound
	\begin{align}\label{eq:reg_total_err_wendland}
		\big\| \big(\calK - \wh \calK_\lambda^\mathrm{MC}\big) f \big\|_\infty \leq C_1 \left( h_\calX^{l+1/2} + \lambda^{1/2} \right) + \veps\cdot C_2\cdot \sqrt{\frac{d}{\lambda}}
	\end{align}
	for every $\veps > 0$ and all $f \in H^{\sigma_{n,l}}(\Omega)$ with $\|f\|_{H^{\sigma_{n,l}}} \leq 1$ with probability at least $1 - 2 d\exp\left( -\frac{m\varepsilon^2}{2 C_2^2}\right)$ and constants $C_1, C_2 > 0$ only depending on $\kernel_{n,l}$.
\end{corollary}
\begin{proof}
	The proof of the first claim is analogous to the proof of \Cref{thm:total_error_bound} by splitting 
	\begin{align*}
		\| (\calK - \wh \calK_\lambda^\mathrm{MC}) f \|_\infty \leq \| (\calK - \wh \calK_\lambda) f \|_\infty + \| (\wh\calK_\lambda - \wh \calK_\lambda^\mathrm{MC}) f \|_\infty
	\end{align*}
	and bounding the deterministic first part by~\Cref{prop:err_tikh}. The estimation error in the second part is bounded by~\Cref{cor:mc_error1_tikh} with the respective probability and yields \eqref{eq:reg_total_err}. 
	
	The second claim combines the second part of \Cref{prop:err_tikh} for the deterministic error bound with setting $C_2 = \|\kernel_{n,l}\|_\infty^{1/2}$, obtaining \eqref{eq:reg_total_err_wendland}. 
\end{proof}
The Matérn kernel analogue to the Wendland special case of \Cref{cor:reg_total_err_bound} is stated in \Cref{cor:reg_total_err_bound_matern} in \Cref{app:matern_kernels}.

\section{Numerical examples}\label{ch:numerics}
In this part, we verify the kernel EDMD method for stochastic systems $Y \mid \{X = x\}$ by considering numerical examples from stochastic differential equations (SDEs). For this, we assume that \Cref{ass:rkhs_invariance} holds and that we can sample from the conditional distribution $\rho_{x}$ of $Y$ given $X = x$. We simulate from SDEs of the form
\begin{equation} \label{e:sde_chap5}
	\begin{aligned}
		d X_t &= b(X_t) \,dt\, + \sigma(X_t) \,d W_t, \\
		X_0 &= x_0 \in \Omega,
	\end{aligned}  
\end{equation}
with Lipschitz continuous drift $b:\R \to \R$ and diffusion $\sigma: \R \to \R$, and Brownian motion $(W_t)_{t \geq 0}$, via the Euler-Maruyama method (see, e.g., \cite{KloePlat1992}). \cite[Theorem 5.2.1]{Okse2000} guarantees the existence of a unique solution $(X_t)_{t \geq 0}$. A basic example for this case is the \emph{Ornstein-Uhlenbeck process} (OU process) given by the linear SDE
\begin{equation}\label{e:ou_sde}
	\begin{aligned}
		d X_t &= -\theta X_t \,dt\, + \sigma \,d W_t, \\
		X_0 &= x_0 \in \Omega,
	\end{aligned}
\end{equation}
with drift coefficient $\theta > 0$ and diffusion coefficient $\sigma > 0$.

We fix a timestep $T > 0$, set $X = X_0$ and $Y = X_T$. To apply the setting of \Cref{ch:hs_invar,ch:err_bounds}, the solution process should remain in some bounded domain~$\Omega = (a,b)$ until time $T$. \cite{PhilSchaWort2025} elaborates on a sufficient condition on the drift and diffusion for the invariance of a compact interval under the SDE flow. Another way to guarantee the boundedness is by reflection at the boundaries. For this, we consider the \emph{doubly reflected} OU process on $\Omega = (0,1) $. It is formally given by the SDE
\begin{equation}\label{e:sde_drou_chap5}
	\begin{aligned}
		d X_t &= -\theta X_t \,dt + \sigma \, dW_t + d L_t - dU_t, \\
		X_0 &= x_0 \in (0,1),
	\end{aligned}
\end{equation}
where $L_t \geq 0$ and $U_t \geq 0$ are the minimal nondecreasing processes that ensure $X_t \geq 0$ and $X_t \leq 1$ with reflection at the boundaries for~$t \geq 0$, respectively. The processes $L$ and $U$ can be interpreted as~``idleness process'' and ``loss process''; they satisfy $\int_0^T \chi_{\{X_t > 0\}} \,d L_t = 0$ and $\int_0^T \chi_{\{X_t < 1\}} \,d U_t = 0$ for any $T > 0$, respectively. We refer to \cite{HuanMandSpre2014} for details. 

As an example of a nonlinear stochastic differential equation, we consider the \emph{double-well potential} $V(x) = \frac{1}{2}x^2 - \frac{1}{4} x^4$ and set the drift to $b(x) = V'(x)$, yielding the SDE
\begin{equation}\label{eq:sde_dw_chap5}
	\begin{aligned}
		d X_t &= (X_t - X_t^3) \,dt\, + \sigma \,d W_t, \\
		X_0 &= x_0 \in \Omega,
	\end{aligned}
\end{equation}
with $\sigma > 0$. To this process, we apply an affine transformation to encapsulate the behavior at the wells $\pm 1$ inside the boundaries $a = 0$ and $b = 1$. We also reflect at those boundaries as described above, thereby bounding the derivative of the drift function, ensuring its Lipschitz continuity. 

\Cref{alg:stochastic_kedmd_theory} summarizes our method to approximate the action of the stochastic Koopman operator on some observable $f: \Omega \to \bR$.

\begin{algorithm}
	\caption{Kernel EDMD for Stochastic Koopman Operator Approximation}
	\label{alg:stochastic_kedmd_theory}
	\begin{algorithmic}[1]
		
		\Require
		Design points $\mathcal{X}=\{x_1,\dots,x_d\}$;
		kernel $\kernel$;
		observable $f$;
		training parameter $m_{\mathrm{train}}$;
		regularization parameter $\lambda$;
		evaluation grid $(x_\ell^{\mathrm{grid}})_{\ell=1}^L$.
		
		\Statex \textbf{Data matrices}
		\State Form the kernel matrix
		\[
		\mathbf{K}_{\mathcal{X}} = \bigl(\kernel(x_i,x_j)\bigr)_{i,j=1}^d.
		\]
		\State Define $f_{\mathcal{X}} = (f(x_1),\dots,f(x_d))^\top$.
		
		\Statex \textbf{Sampling step}
		\State For each $j=1,\dots,d$ draw i.i.d.\ samples 
		\[
		y_j^{(m)} \sim \rho_{x_j},
		\qquad m=1,\dots,m_{\mathrm{train}}.
		\]
		
		\State Estimate the propagated kernel matrix
		\[
		(\bbK_{\calX, \calY}^\mathrm{MC})_{ij}
		= \frac{1}{m_{\mathrm{train}}} 
		\sum_{m=1}^{m_{\mathrm{train}}}
		\kernel(x_i, y_j^{(m)}).
		\]
		
		\Statex \textbf{Koopman estimate}
		\For{$\ell = 1,\dots,L$}
		\State Form kernel vector 
		\[
		\bbk_{\mathcal{X}}(x_\ell^{\mathrm{grid}})
		= (\kernel(x_1,x_\ell^{\mathrm{grid}}),\dots,\kernel(x_d,x_\ell^{\mathrm{grid}}))^\top.
		\]
		\State Compute the estimated Koopman action
		\[
		\widehat{\mathcal{K}}^\mathrm{MC} f(x_\ell^{\mathrm{grid}})
		=
		f_{\mathcal{X}}^\top\,
		\bbK_\calX^{-1}\,
		\bbK_{\calX, \calY}^\mathrm{MC}\, 
		\bbK_\calX^{-1}\,
		\bbk_{\mathcal{X}}(x_\ell^{\mathrm{grid}}).
		\]
		\EndFor
		
		\State \Return $\widehat{\mathcal{K}}^\mathrm{MC}f(x_\ell^{\mathrm{grid}})$ for $\ell=1,\dots,L$.
		
	\end{algorithmic}
\end{algorithm}

In the setting of the stochastic Koopman operator \eqref{eq:koop_identity}, we set $X = X_0$ and $Y = X_{T}$. This means that we are interested in the evolution of observables $f: \Omega \to \bR$ under the stochastic process that solves \eqref{e:sde_chap5}, i.e., 
\begin{align}\label{eq:koop_identity_sde}
	\calK f(x) = \bE[f(X_{T}) \mid X_0 = x].
\end{align}
For numerical testing, we can generally proceed as usual in kernel EDMD methods using the approximation \eqref{eq:mcedmd}, that is,
\begin{align*}
	\widehat \calK^{\mathrm{MC}}f = f_\calX^\top \mathbf{K}_{\calX}^{-1}\bbK_{\calX, \calY}^{\mathrm{MC}} \mathbf{K}_{\calX}^{-1}\bbk_\calX,
\end{align*}
or the regularized approximant \eqref{eq:regularized_mcemdmd}, i.e.,
\begin{align*}
	\widehat \calK^{\mathrm{MC}}f = f_\calX^\top (\mathbf{K}_{\calX} + \lambda I)^{-1}\bbK_{\calX, \calY}^{\mathrm{MC}} (\mathbf{K}_{\calX} + \lambda I)^{-1}\bbk_\calX.
\end{align*}
Aside from the observable $f: \Omega \to \bR$ and design points $\calX = \{x_1, \dots, x_d\} \subset \Omega$ that need to be chosen beforehand, we introduce the parameter $m_\mathrm{train} \in \bN$ that specifies the number of samples we draw from the conditional distributions $\rho_{x_j}$, $j= 1, \dots, d$. 
We note that the error analysis in the stochastic kEDMD setting presents several challenges compared to the deterministic setting:
\begin{enumerate}
	\item Conditional expectations of the form $\calK f(x_0) = \bE[f(X_{T}) \mid X_0 = x_0]$ for some $x_0 \in \Omega$ do not necessarily admit an analytical expression. Therefore, we usually do not have access to the exact action of the Koopman operator on the observable $f$ ("ground truth"). Instead, we use a further Monte Carlo approximation 
	\begin{align}\label{eq:mc_testing}
		\calK f(x_\mathrm{test}) \approx \frac{1}{m_\mathrm{truth}} \sum_{l = 1}^{m_\mathrm{test}} f(y_l)
	\end{align}
	where $y_l \overset{\text{i.i.d.}}{\sim} \rho_{x_\mathrm{test}}$ are i.i.d.\ samples of the simulated SDE solution at time $T$ and $0 \ll m_\mathrm{truth} \in \bN$ is our ``ground truth parameter''. It should also be noted that we introduce a discretization error by simulating the SDE via the Euler-Maruyama method. To control this error, a fine enough grid on $\Omega = (0,1)$ is employed. 
	\item Even for an analytically calculated ground truth, we can only approximate the $L^\infty$-error over some grid on $\{x_1^\mathrm{grid}, \dots, x_L^\mathrm{grid} \} \subset \Omega$. This gives us an ``empirical'' approximated $L^\infty$-error $\wh{\mathrm{err}} = \max_l |\calK f(x_l^\mathrm{grid}) - \wh{\calK}^\mathrm{MC} f(x_l^\mathrm{grid})|$. 
	\item The randomness of the propagated matrix $K_{\calX, \calY}^\mathrm{MC}$ necessitates a comparison over multiple realizations of the Koopman approximation. In our SDE setting, this means that we build matrices $\bbK_{\calX, \calY}^\mathrm{MC, n}$, $n = 1, \dots, n_\mathrm{pred}$, and use the mean values and standard deviations of the empirical $L^\infty$-errors $\wh{\mathrm{err}}^{(n)}$ per realization to judge the performance of the approximation $\wh \calK^\mathrm{MC}$. 
\end{enumerate}

\Cref{thm:total_error_bound} decomposes the $L^\infty$-error of the stochastic Koopman approximation into a deterministic part controlled by the fill distance $h_\calX$ and a probabilistic part influenced by the Monte Carlo sample size, which we denote by $m_\mathrm{train}$ to differentiate it from the number of Monte Carlo samples $m_\mathrm{truth}$ used for ground truth approximation. 
Throughout the testing procedure, we frequently sample the SDE solution with initial values $x_0 \in \Omega$ at time $T > 0$ from the conditional distribution $\rho_{x_0, T}$. First, we start trajectories at chosen evaluation grid points $\calX^\mathrm{grid} \subset \Omega$ for ground truth computation. Later, we repeatedly simulate trajectories starting at the design points $\calX$ to construct realizations $\bbK_{\calX, \calY}^\mathrm{MC}$ in the Koopman approximation $\wh \calK^\mathrm{MC}$. For numerical demonstration on the domain~$\Omega = (0,1)$, we first chose fill distances $h_\calX$ such that $d = 2/h_\mathcal{X} \in \bN$. The fill distance was then achieved by spacing out the $d$ design points $\calX = \{x_1, \dots, x_d\}$ evenly on $\Omega = (0,1)$: We set $x_1 = h_\mathcal{X}$, $x_d = 1 - h_\mathcal{X}$, and $x_i = (2i - 1)h_\mathcal{X}$ for~$i = 2, \dots, d-1$.
\Cref{alg:stochastic_kedmd_numerics} illustrates the procedure.

\begin{algorithm}[htb]
	\caption{Error Analysis for the Stochastic Kernel EDMD Koopman Approximation}
	\label{alg:stochastic_kedmd_numerics}
	\begin{algorithmic}[1]
		
		\Require \hfill
		\begin{itemize}
			\item number of Koopman approximation realizations $n_{\mathrm{pred}} \in \mathbb{N}$;
			\item skEDMD parameters = \{fill distance $h_\calX$; kernel $\kernel$; training parameter $m_{\mathrm{train}}$; observable $f$; regularization parameter $\lambda$; evaluation grid $\calX^\mathrm{grid} = \{x_\ell^{\mathrm{grid}}\}_{\ell=1}^L$\};
			\item simulation parameters = \{SDE parameters; SDE simulator; timestep T\};
			\item ground truth parameters.
		\end{itemize}
		
		\Statex \textbf{Ground truth}
		\For{$\ell=1,\dots,L$}
		\State For $m = 1, \dots, m_\mathrm{truth}$, simulate $y_\ell^{(m)} \sim \rho_{x_\ell, T}$ with simulation parameters and set
		\[ \mathcal{K}_{\mathrm{truth}} f(x_\ell^{\mathrm{grid}}) = \frac{1}{m_{\mathrm{truth}}} \sum_{m=1}^{m_{\mathrm{truth}}} y_\ell^{(m)}
		\]
		\EndFor
		
		\Statex \textbf{Design points}
		\State Generate equidistant $\calX = \{x_i\}_{i=1}^d$ for fill distance (at least) $h_{\mathcal{X}}$.
		
		\Statex \textbf{Kernel EDMD iterations}
		\For{$n = 1,\dots,n_{\mathrm{pred}}$}
		
		\State Apply \Cref{alg:stochastic_kedmd_theory} with skEDMD and simulation parameters
		to obtain approximations $\widehat{\mathcal{K}}^{\mathrm{MC}} f(x_\ell^{\mathrm{grid}})$ for $l = 1, \dots, L$
		
		\State Compute approximation of the $L^\infty$-error over the evaluation grid: 
		\[
		\widehat{\mathrm{err}}_n = \max_\ell \, \lvert \mathcal{K}_{\mathrm{truth}} f(x_\ell^{\mathrm{grid}})
		- \widehat{\mathcal{K}}^{\mathrm{MC}} f(x_\ell^{\mathrm{grid}})\rvert
		\]
		\EndFor
		
		\State \Return mean and standard deviation \ of $\widehat{\mathrm{err}}_n$, $n=1,\dots,n_{\mathrm{pred}}$.
		
	\end{algorithmic}
\end{algorithm}

For numerical stability, we apply a Tikhonov regularization with parameter $\lambda$ as in \eqref{eq:regularized_mcemdmd} to alleviate the bad conditioning of kernel matrices $\bbK_\mathcal{X}$. The parameters $m_\mathrm{train}$ and $\lambda$ are occasionally adapted to the chosen fill distance $h_\calX$ as indicated in the table of \Cref{fig:linf_ou_const_vs_adapt_left_right}. 

To settle on a number of realizations $n_\mathrm{pred}$ to control for the randomness in~$\wh  \calK^\mathrm{MC}$, we conducted numerical tests for different choices of $h_\calX$ and $m_\mathrm{train}$. Our tests showed that the standard deviation of the empirical $L^\infty$-error seems to stabilize after as little as $n_\mathrm{pred} = 30$ realizations. Therefore, we average approximation errors over $30$ realizations in the following. 

The remainder of this section studies the error behavior for the linear Ornstein-Uhlenbeck process and the nonlinear process governed by the double-well potential in dependence of the fill distance $h_\calX$ and Monte Carlo sample size $m_\mathrm{train}$. In all cases, we consider the nonlinear observable $f(x) = e^x$. We chose to compare the
\begin{align*}
	\text{Wendland kernel}\quad \kernel_{1, 1}^\text{wendland}(x,y) &= (1 - r)_+^3 \cdot (3r + 1),~\text{and} \\
	\text{Mat\'{e}rn kernel}\quad \kernel_{3/2}^\text{mat\'{e}rn}(x,y) &= \left(1+ \sqrt{3}r\right) \exp \left( - \sqrt{3}r \right),
\end{align*} 
where we denote $r = \|x-y\|_2$.

\subsection{Linear SDE: Ornstein-Uhlenbeck process}\label{subch:numerics_ou}
In the linear Ornstein-Uhlenbeck SDE \eqref{e:ou_sde}, we set $\theta = 1$, $\sigma = 0.2$ and $T = 1$ to simulate from
\begin{equation}\label{eq:ou_sde_doubly_reflected}
	\begin{aligned}
		d X_t &= - X_t \,dt\, + 0.2 \,d W_t + d L_t - d U_t, \\
		X_0 &= x_0 \in \Omega,
	\end{aligned}
\end{equation}
with initial value $x_0$ and auxiliary processes $L$ and $U$ to achieve reflection at the boundaries of $\Omega = (0,1)$. 

As laid out in \Cref{alg:stochastic_kedmd_numerics}, we approximate the $L^\infty$-error per realization as maximal absolute error over all evaluation grid points. We then average over all $n_\mathrm{pred}$ realizations to obtain an empirical $L^\infty$-error and compute its standard deviation.

\begin{figure}[htb]
	\centering
	
	\begin{minipage}{0.68\textwidth}
		\centering
		\includegraphics[width=\linewidth]{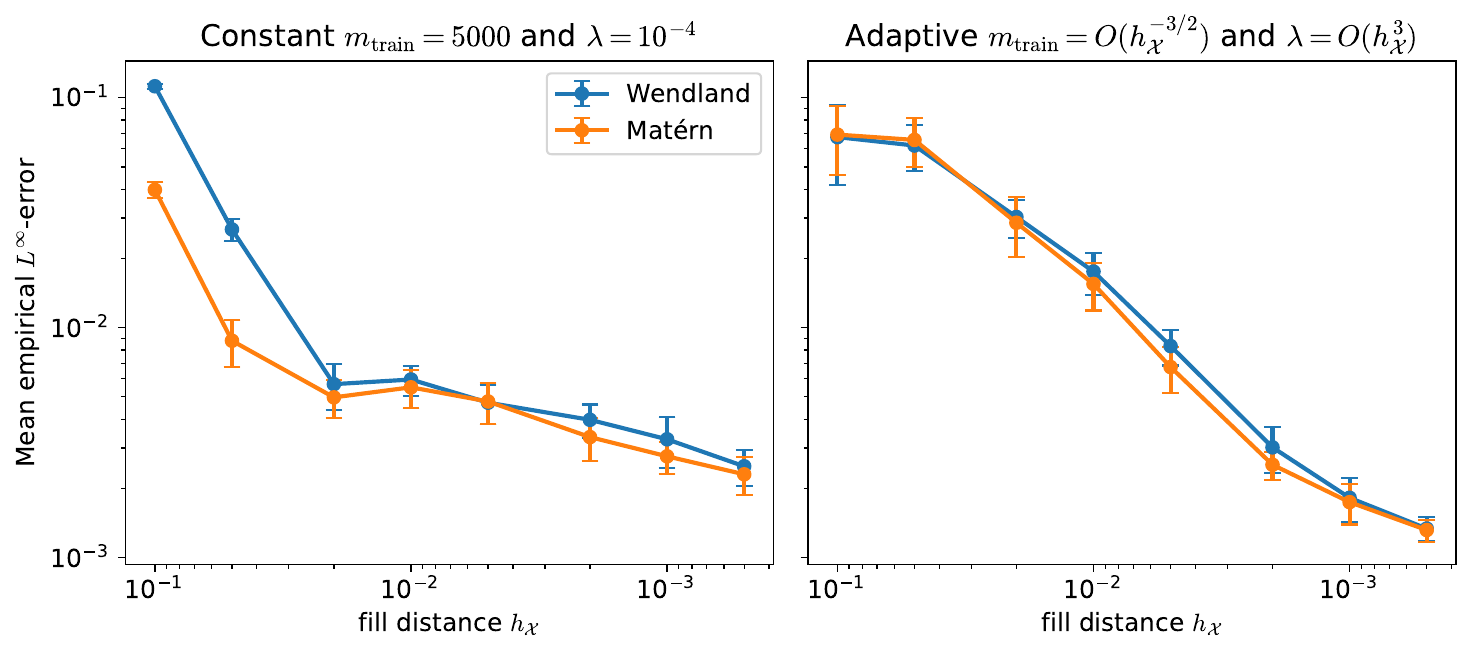}
	\end{minipage}
	\hfill
	\begin{minipage}{0.27\textwidth}
		\vspace{-0.8\baselineskip}
		\centering
		\smallskip
		\footnotesize
		\textbf{Adaptive parameters}
		\setlength{\tabcolsep}{4pt}
		
		\begin{tabular}{c c c}
			\toprule
			\boldmath$h_\mathcal{X}$ & \boldmath$\lambda$ & \boldmath$m_\mathrm{train}$ \\
			\midrule
			$1 \cdot 10^{-1}$  & $1\cdot 10^{-3}$  & 20    \\
			$5 \cdot 10^{-2}$  & $\frac{5}{4}\cdot 10^{-4}$  & 30    \\
			$2 \cdot 10^{-2}$  & $8\cdot 10^{-6}$  & 127   \\
			$1 \cdot 10^{-2}$  & $1\cdot 10^{-6}$  & 387   \\
			$5 \cdot 10^{-3}$  & $\frac{5}{4}\cdot 10^{-7}$  & 1183  \\
			$2 \cdot 10^{-3}$  & $8\cdot 10^{-9}$  & 5166  \\
			$1 \cdot 10^{-3}$  & $1\cdot 10^{-9}$  & 15684 \\
			$5 \cdot 10^{-4}$  & $\frac{5}{4}\cdot 10^{-10}$ & 47428 \\
			\bottomrule
		\end{tabular}
	\end{minipage}
	\caption{\footnotesize
		Ornstein--Uhlenbeck process: mean empirical $L^\infty$-error of $n_\mathrm{pred}=30$ realizations of the stochastic Koopman approximation as a function of the fill distance $h_\mathcal{X}$. The plot compares constant and adaptive choices of $m_\mathrm{train}$ and $\lambda$; the table reports the adaptive parameter values used for each $h_\mathcal{X}$. Error bars indicate standard deviation over the realizations.
	}
	\label{fig:linf_ou_const_vs_adapt_left_right}
\end{figure}

The left panel in \Cref{fig:linf_ou_const_vs_adapt_left_right} is a first display of the dependence of the deterministic part of the error on the fill distance for both Matérn and Wendland kernel. However, after an initially steeper decline of error with fill distance, the curve seems to flatten where we would expect a linear decay in the log-log-plot from the deterministic error bound of \Cref{thm:err1}.

\begin{figure}[htb]
	\centering
	\includegraphics[width=.8\linewidth]{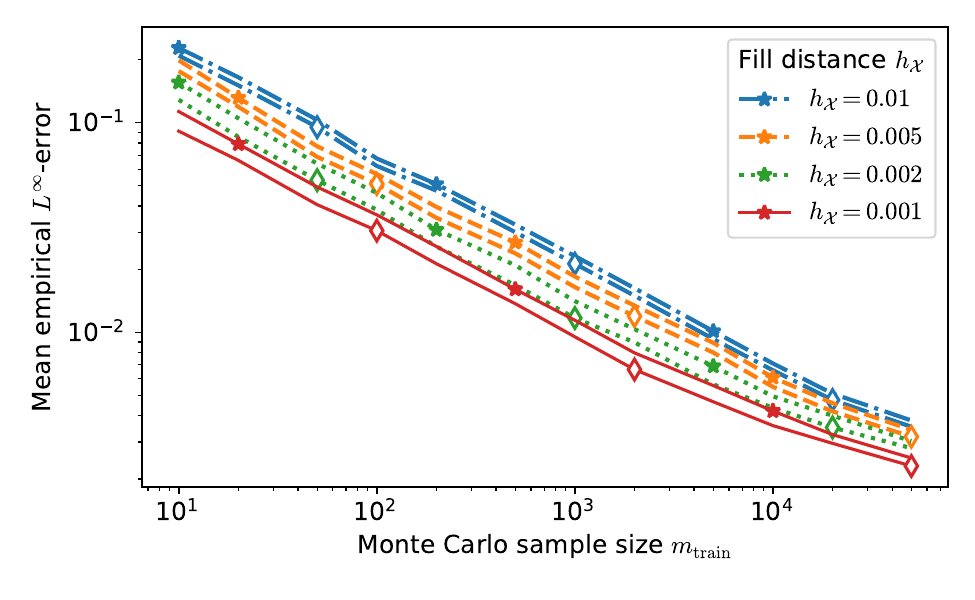}
	\caption{
		\footnotesize
		Ornstein-Uhlenbeck process: mean empirical $L^\infty$-error of $n_\mathrm{pred} = 30$ realizations of the stochastic Koopman approximation in dependence of the number of Monte Carlo simulations $m_\mathrm{train}$. $\star$ and $\diamond$ indicate Wendland and Matérn kernels, respectively.}
	\label{fig:wendland_vs_matern_mtrain_4fd_ou}
\end{figure}

In order to guarantee decay with the fill distance with the suggested rate, we have to adapt $m_\mathrm{train}$ to the fill distance to ensure the probabilistic error behaves similarly and does not become dominant for smaller fill distances. \Cref{app:numerical_considerations} expands on that balancing. We heuristically chose $m_\mathrm{train}(h_\calX) = O(h_\calX^{-3/2}) = O(d^{3/2})$ and $\lambda(h_\calX) = O(h_\calX^{3}) = O(d^{-3})$. 
The right panel of \Cref{fig:linf_ou_const_vs_adapt_left_right} depicts how decreasing fill distances combined with adapted Monte Carlo parameters $m_\mathrm{train}$ and regularization parameters $\lambda$ result in smaller empirical $L^\infty$-errors. The chosen growth rate of the Monte Carlo parameter $m_\mathrm{train}$ in dependence of the fill distance $h_\calX$ is far more conservative than the theoretically derived one in \Cref{app:numerical_considerations}. Still, after adapting the parameters $m_\mathrm{train}$ and $\lambda$, the decay appears close to linear on the log-scale. The Wendland kernel and Matérn kernel perform similarly well and both achieve a comparable decay of the $L^\infty$-error. 

Another insight is gained by fixing a fill distance and only increasing $m_\mathrm{train}$, decreasing the probabilistic error with fixed deterministic error bound. The probabilistic effect is illustrated in \Cref{fig:wendland_vs_matern_mtrain_4fd_ou} for various fill distances. For a fixed $m_\mathrm{train}$, a smaller fill distance results in a smaller $L^\infty$-error. This is consistent with the error bounds on the deterministic, probabilistic and full $L^\infty$-error proven above. Again, the performance of Wendland and Matérn kernel is comparable.

\subsection{Nonlinear SDE: Rescaled double-well potential}\label{subch:numerics_dw}
For the nonlinear SDE with double\=/well potential, we also set $\sigma = 0.2$ and $T = 1$. For better comparability, especially in terms of the fill distance, we want to restrict the corresponding SDE to the same domain~$\Omega = (0,1)$. Since the wells in the regular double-well SDE \eqref{eq:sde_dw_chap5} are at $\pm 1$, we apply the affine transformation $Z_t = (X_t + 2)/4$ to mimic the behavior of the original double-well potential on $(-2, 2)$ on $\Omega = (0,1)$. Solving for $X_t$ gives $X_t = 4 Z_t - 2$, so that the chain rule yields the SDE
\begin{equation}\label{eq:doublewell_sde_doubly_reflected}
	\begin{aligned}
		d Z_t &= \frac{1}{4} \left[(4 Z_t - 2) - (4 Z_t -2)^3 \right] \,dt\, + \frac{\sigma}{4} \,d W_t + d L_t - d U_t, \\
		Z_0 &= z_0 \in \Omega,
	\end{aligned}
\end{equation}
with the auxiliary processes $L$ and $U$ achieving reflection at the boundaries $a = 0$ and $b = 1$. 

\begin{figure}[htb]
	\centering
	\includegraphics[width=1\linewidth]{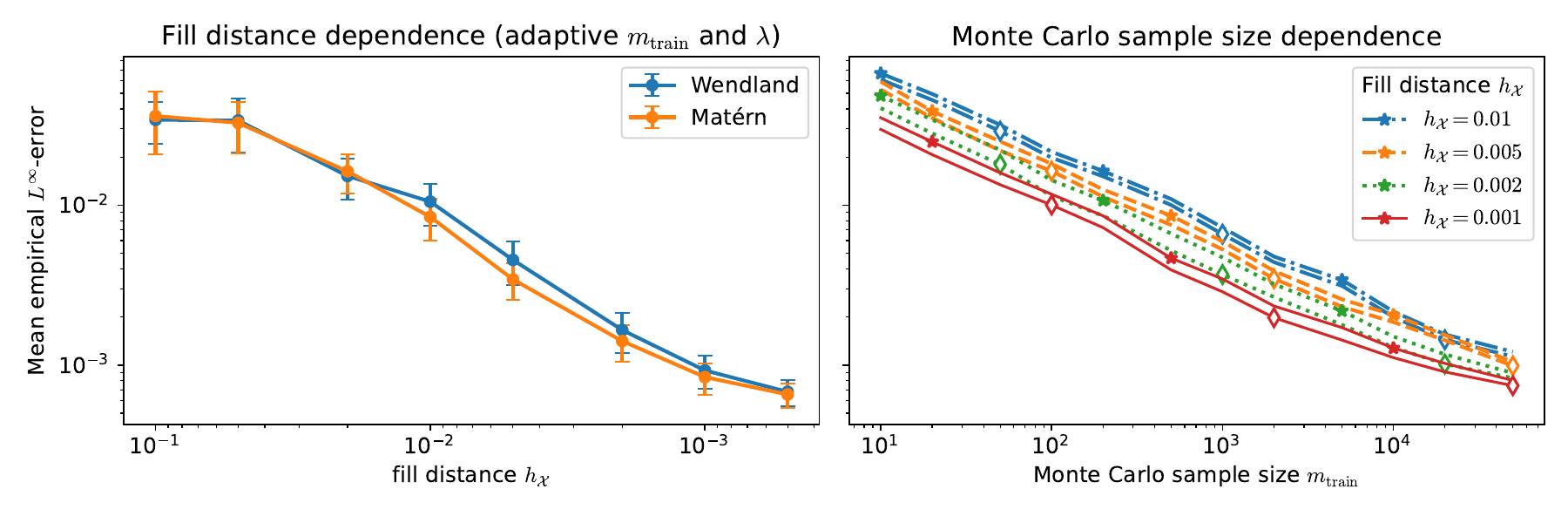}
	\caption{
		\footnotesize
		Double-well SDE: mean empirical $L^\infty$-error with standard deviation of $n_\mathrm{pred} = 30$ realizations of the stochastic Koopman approximation in dependence of fill distance $h_\calX$ and Monte Carlo sample size $m_\mathrm{train}$. Parameters in the left panel were adapted as indicated by the table in \Cref{fig:linf_ou_const_vs_adapt_left_right}. In the right panel, $\star$ and $\diamond$ indicate Wendland and Matérn kernels, respectively.}
	\label{fig:linf_vs_fd_mtrain_4fd_dw}
\end{figure}

As \Cref{fig:linf_vs_fd_mtrain_4fd_dw} shows, the double-well SDE demonstrates similar error behavior. Mirroring the adaptive parameter regime from the linear OU process above, we scaled $m_\mathrm{train}$ and $\lambda$ with the fill distance to demonstrate the decay of the $L^\infty$-error. The decrease of the probabilistic error for various fixed fill distances also follows a similar pattern as in \Cref{fig:wendland_vs_matern_mtrain_4fd_ou}. 
Furthermore, we see that, while comparable in performance, Matérn kernels seem to have a slight edge over Wendland kernels for each given fill distance.

\section{Conclusions and Outlook}\label{ch:conclusion}
\noindent We proposed a Monte Carlo-based estimator in a kernel EDMD approach for the Koopman operator of stochastic systems and analyzed its pointwise error. Extending the results of \cite{KohnPhil24}, we decompose the error into the well-known deterministic projection error and the probabilistic error introduced by the entry-wise estimation of expected values in the exact ``propagation matrix'' $\bbK_{\calX, \calY}$. To apply known projection error bounds, we derive sufficient conditions for the invariance of (fractional) Sobolev spaces $H^s(\Omega)$ under the stochastic Koopman operator, which are induced as RKHS by Wendland and Matérn kernels. Novel probabilistic error bounds for our estimator are proven to obtain a pointwise bound on the full approximation error with high probability. 

In this work, we only considered two exemplary kernels, i.e., Mátern and Wendland. In future work, the respective error bounds may be derived for other kernels and their corresponding RKHSs. 
Furthermore, a canonical (and important) extension is the one to \textit{multilevel Monte Carlo} (MLMC) to increase sample and, thus, computational efficiency by variance reduction.

\bibliographystyle{plain}
\bibliography{References}

\appendix

\section{Proof of Monte Carlo error bound from \texorpdfstring{\Cref{rem:l2_l2_mc_error_bound}}{Remark 4.7}}\label{app:l2_l2_MC}
Applying Hölder's inequality according to \Cref{rem:l2_l2_mc_error_bound} with $p = q = 2$ instead of the approach $p = \infty$, $q = 1$ from \Cref{thm:mc_error1}, one obtains
\begin{align*}
	|(\widehat \calK^{\mathrm{MC}}-\widehat \calK)f(x)| \leq \|g-g^{\mathrm{MC}}\|_2 \|K_\calX^{-1} \mathbf{k}_\calX(x)\|_2.
\end{align*}
We then have
\begin{align*}
	\bP\big(\|g-g^\mathrm{MC}\|_2\ge\veps\|f\|_H)
	&= \bP\left(\sum_{j=1}^d\left|\frac{1}{m_j}\sum_{l=1}^{m_j} \big[(S_\calX f)(y_{l,j}) - \bE_j\big]\right|^2\,\ge\,\veps^2\|f\|_H^2\right)\\
	&\le \bP\left(\exists j\in \{1,\dots,d\} : \left|\frac{1}{m_j}\sum_{l=1}^{m_j} \big[(S_\calX f)(y_{l,j}) - \bE_j\big]\right|^2\,\ge\,\frac{\veps^2\|f\|_H^2}{d}\right)\\
	&\le \sum_{j=1}^d \bP\left(\left|\sum_{l=1}^{m_j} \big[(S_\calX f)(y_{l,j}) - \bE_j\big]\right|\,\ge\,\frac{\veps m_j\|f\|_H}{\sqrt d}\right)\\
	&\le 2\sum_{j=1}^d\exp\left(-\frac{m_j\veps^2}{2d\|\kernel\|_\infty}\right).
\end{align*}
Hence, if $m_j=m$ for all $j\in \{1,\dots,d\}$, then
\[
\bP\big(\|g-g^\mathrm{MC}\|_2\ge\veps\|f\|_H)\,\le\,2d\exp\left(-\frac{m\veps^2}{2d\|\kernel\|_\infty}\right).
\]
We may improve this by using the Hilbert-space Hoeffding inequality. For this, let $m_j = m$ for all $j$, and set $\xi_l = [(S_\calX f)(y_{l,j})]_{j=1}^d - g\in\R^d$, $l\in \{1, \dots, m\}$. Then, we have $\bE[\xi_l] = 0$ and 
\[
\|\xi_l\|_2^2 = \sum_{j=1}^d \left|(S_\calX f)(y_{l,j}) - \int (S_\calX f)(y)\,\rho(x,dy)\right|^2.
\]
Since $|(S_\calX f)(y)| \leq \|\kernel\|_\infty^{1/2}\|f\|_H$, we obtain $\|\xi_l\|_2^2\le 4d\|\kernel\|_\infty\|f\|_H^2$ for all $l\in \{1,\dots, m\}$. Also,
\begin{align*}
	\bP\big(\|g-g^\mathrm{MC}\|_2\ge\veps\|f\|_H)
	&= \bP\left(\sum_{j=1}^d\left|\frac{1}{m}\sum_{l=1}^{m} \xi_{l,j}\right|^2\,\ge\,\veps^2\|f\|_H^2\right) = \bP\left(\left\|\frac{1}{m}\sum_{l=1}^{m}\xi_l\right\|\,\ge\,\veps\|f\|_H\right)\\
	&\le 2\exp\left(-\frac{m\veps^2}{8d\|\kernel\|_\infty}\right).
\end{align*}
We are left with estimating $\|\bbK_\calX^{-1}\bbk_\calX(x)\|_2$. The map $\tau : V_\calX\to\R^d$, $\tau f = f_\calX$, is surjective (namely, if $v\in\R^d$ is given, then $\tau(v^\top\bbK_\calX^{-1}\bbk_\calX) = v$), hence
\begin{align*}
	\|\bbK_\calX^{-1}\bbk_\calX(x)\|_2
	&= \sup_{\|v\|_2=1}|v^\top\bbK_\calX^{-1}\bbk_\calX(x)| = \sup_{f\in V_\calX,\,f\neq 0}\frac{|f_\calX^\top\bbK_\calX^{-1}\bbk_\calX(x)|}{\|f_\calX\|_2} = \sup_{f\in V_\calX,\,f\neq 0}\frac{|f(x)|}{\|f_\calX\|_2}\\
	&= \sup_{f\in V_\calX,\,f\neq 0}\frac{|\<f,\Phi_x\>_H|}{\|f_\calX\|_2}\,\le\,\|\kernel\|_\infty^{1/2}\sup_{f\in V_\calX,\,f\neq 0}\frac{\|f\|_H}{\|f_\calX\|_2} = \|\kernel\|_\infty^{1/2}\|\bbK_\calX^{-1}\|_{2\to 2}^{1/2}.
\end{align*}
Alternatively, by the $\ell^1$-$\ell^2$-inequality, one may bound the deterministic term as in \Cref{thm:mc_error1} by $\|\bbK_\calX^{-1}\bbk_\calX(x)\|_2 \leq \|\bbK_\calX^{-1}\bbk_\calX(x)\|_1 \leq \sqrt{\max_{v\in\one}v^\top \mathbf{K}_\calX^{-1}v}$.

\section{Recap: Wendland kernels}\label{app:wendland_kernels}
Wendland kernels are constructed using particular compactly supported \emph{radial basis functions}, which are constructed in the following way:
For $m \in \bN$, define the function $\beta_m: \bR_+ \to \bR, \beta_m(r) = (1 - r)_+^m {= \max(1-r, 0)^m}$. Given the operator $\calI$ that acts on measurable non-negative functions $\beta: \bR_+ \to \bR_+$ by 
\begin{align*}
	\calI \beta(r) = \int_r^\infty t \beta(t)\,dt,
\end{align*}
and given dimension $n \geq 1$, we choose a smoothness $l \in \bN$ and set $\beta_{n, l} = \calI^l \beta_{\lfloor \frac{n}{2}\rfloor + l + 1}$. This function is of the form
\begin{align*}
	\beta_{n,l}(r) = \begin{cases}
		p_{n,l}(r),&0 \leq r \leq 1, \\
		0, &r > 1
	\end{cases}
\end{align*}
with a univariate polynomial $p_{n,l}$ of degree $\lfloor \frac{n}{2} \rfloor + 3l + 1$ whose coefficients can be computed recursively by~\cite[Theorem 9.12]{Wend04}. We can now define the \emph{Wendland radial basis function} $B_{n,m}: \bR^n \to \bR$ of dimension $n$ and smoothness degree $m$ by
\begin{align*}
	B_{n, l}(x) := \beta_{n,l}(\|x\|_2) 
\end{align*}
and the corresponding \emph{Wendland kernel} by $\kernel_{n,l}: \bR^n \times \bR^n \to \R, \kernel_{n,l}(x, y) = B_{n,l}(x-y) = \beta_{n,l}(\|x-y\|_2)$. 

\section{Error bounds for Matérn kernels}\label{app:matern_kernels}
Similarly to the Wendland kernels, Matérn kernels are constructed via radial basis functions. Specifically, the \emph{Matérn covariance function} $\calM_{\nu, \alpha}: \bR_+ \to \bR$ with smoothness $\nu > 0$ and length-scale parameter $\alpha > 0$ is given by
\begin{align*}
	\calM_{\nu, \alpha}(r) =   \frac{2^{1 - \nu}}{\Gamma(\nu)} \left( \frac{\sqrt{2\nu}\, r}{\alpha} \right)^\nu \kernel_\nu\left(\frac{\sqrt{2\nu}\, r}{\alpha} \right),
\end{align*}
where $\Gamma$ denotes the gamma function and $\kernel_\nu$ denotes the modified Bessel function of the second kind.
Denoting the Euclidean distance of $x, y \in \bR^n$ by $\| \cdot\|_2$, the Matérn kernel $\kernel_{\nu, \alpha}: \bR^n \times \bR^n \to \R$ is given by
\begin{align*}
	\kernel_{\nu, \alpha}(x,y) = \calM_{\nu, \alpha}(\|x-y\|_2).
\end{align*}
The parameter $\nu$ controls the differentiability of the Matérn kernel at $x = y$ through $\kernel_{\nu, \alpha} \in C^{2 \lfloor \nu \rfloor}(\R^n \times \R^n)$ – with increasing $\nu$, the Matérn kernels become smoother. Equivalently, the sample paths of a Gaussian process with the Matérn kernel as covariance function are $\lfloor \nu\rfloor$-times mean-square differentiable, which provides a flexible way to control smoothness and makes Matérn kernels a popular choice in this setting. 
Another central property is the convergence 
\begin{align*}
	\calM_{\nu, \alpha}(r) \xrightarrow[]{\nu \to \infty} \exp\left( -\frac{r^2}{2\alpha^2} \right)
\end{align*}
to the squared exponential (Gaussian) kernel, uniformly on compact sets. Of particular interest are the cases $\nu = s + 1/2$ for~$s \in \bN$: In that case, the Matérn covariance function is given as the product of a negative exponential and a polynomial of order $s$ \cite{PorcBeviSchabOate2024}.
Up to scaling constants, the Matérn kernel $\kernel_{\nu, \alpha}$ on $\bR^n$ induces the Sobolev space $H^{\nu + n/2}$ with equivalent norms; for the dimensionless Matérn kernel $\kernel_\nu := \kernel_{\nu, \sqrt{2\nu}}$, the RKHS norm even coincides with the standard Sobolev norm of order $\nu + n/2$ up to a constant factor. We refer the reader to \cite{FassYe11, PorcBeviSchabOate2024} for more insights into Matérn kernels. 

Matérn kernels $\kernel_\nu(x,y) = \calM_{\nu, 1}(\|x-y\|_2)$ of smoothness $\nu > 0$ induce the Sobolev spaces $H^{\nu + n/2}$. Denote the native norm of the RKHS by $\|\cdot\|_\bH$. The Matérn kernel is $2\lfloor \nu \rfloor$-times differentiable. Furthermore, (normalized) Matérn kernels are uniformly bounded by $1$. Thus, according to the discussion after \cite[Theorem 11.13]{Wend04}, the bound 
\begin{align*}
	\| (\calI - S_\calX)f \|_\infty \leq C h_\calX^{ \lfloor \nu \rfloor} \| f\|_\bH 
\end{align*}
with a constant $C > 0$ holds for all $f \in \bH$. 

Using the above properties, we obtain analogous results as in \Cref{ch:err_bounds} for kernel EDMD approximations with underlying Matérn kernels. Following the derivations for Wendland kernels, we start with the deterministic error bound characterized by the projection error. 

\begin{corollary}[Deterministic error bound in Matérn native spaces]\label{cor:det_err_matern}
	Let $\nu > 0$. Denote by $\kernel_{\nu}: \bR^n \to \bR^n$ the (dimensionless) Matérn kernel generating the Sobolev space $H^{\nu + n/2}(\Omega)$ with equivalence between the native norm $\| \cdot \|_\bH$ and the classical Sobolev norm $\| \cdot \|_{H^{\nu + n/2}}$. Assume that $\calK \in L(C_b, C_b)$ and that \eqref{e:hl_charac} is satisfied for some $\ell \in \bN$ with $\ell > \nu + n/2$. Then, the Koopman approximation $\wh\calK$ from \eqref{eq:approx1} satisfies
	\begin{align}\label{eq:det_err_matern}
		\|\calK - \wh \calK \|_{\bH \to C_b} \leq C (1+\|\calK\|_{\bH\to \bH})\cdot h_\calX^{\lfloor \nu \rfloor}.
	\end{align}
\end{corollary}
\begin{proof}
	The Matérn kernel satisfies $\kernel_\nu \in C^{\lfloor 2 \nu \rfloor}(\bR^n \times \bR^n)$.
	According to \cite[Theorem 11.13]{Wend04} and the subsequent discussion, we have
	\begin{align*}
		\| \calI - S_\calX \|_{\bH \to C_b} \leq h_\calX^{\lfloor \nu \rfloor}
	\end{align*}
	since the Matérn kernel's derivatives of order $2 \lfloor \nu \rfloor$ are continuous on $\overline{\Omega} \times \overline{\Omega}$, allowing for a uniform bound on those derivatives on that compact set. 
	The assertion \eqref{eq:det_err_matern} then follows analogously as in the proof of \Cref{cor:det_err_wendland}.
\end{proof}

Since the probabilistic error bound for the Matérn case is virtually identical to the Wendland case, merely inserting $\|\kernel_\nu\|_\infty = 1$, we state the analogue to \Cref{thm:total_err_wendland} for the Matérn case below.

\begin{corollary}[$L^\infty$-error bound for Matérn native spaces]\label{cor:total_error_bound_matern}
	Let the assumptions of \Cref{thm:total_error_bound} hold. Let $\nu > 0$ and denote by $\bH = H^{\nu + n/2}(\Omega)$ the reproducing kernel Hilbert space induced by the Matérn kernel $\kernel_{\nu}$. Then, for every $\veps > 0$  and all $f \in \bH$ with $\|f\|_\bH\le 1$, the full error admits the bound
	\begin{equation}\label{eq:total_err_matern}
		\begin{aligned}
			\| (\calK - \calK^{\mathrm{MC}})f\|_{\infty} \leq\, &C\,(1 + \|\calK\|_{\bH\to \bH})\cdot h_\calX^{\lfloor \nu \rfloor}
			+ \veps\cdot\sqrt{\max_{v\in\one}v^\top \mathbf{K}_\calX^{-1}v}
		\end{aligned}
	\end{equation}
	with probability at least $1 - 2d\exp\left(-\frac{m\veps^2}{2}\right)$ and $\one = \{v\in\R^d : |v_i|=1\}$.
\end{corollary}
\begin{proof}
	The Matérn covariance function is monotonically decreasing for~$r > 0$, admitting its maximum at $\calM_{\nu, \alpha}(0) = 1$, so $\|\kernel_\nu\|_\infty = 1$. Similarly to \Cref{thm:total_err_wendland}, the result now follows immediately from \Cref{thm:total_error_bound}. 
\end{proof}

Applying Tikhonov regularization with some $\lambda > 0$, as is usual in practice, gives rise to the error bound
\begin{align}\label{eq:det_error_tikh_matern}
	\| \calK f - \wh \calK_\lambda f \|_\infty \leq C \left( h_\calX^{\lfloor \nu \rfloor} + \lambda^{1/2} \right) \| f \|_\bH
\end{align}
for~$\bH = H^{\nu + n/2}(\Omega)$, following the line of argument in \Cref{prop:err_tikh}. To conclude, we state the analogue to the second part \Cref{cor:reg_total_err_bound} for Matérn native spaces. 

\begin{corollary}[Regularized $L^\infty$-error bounds for Matérn native spaces]\label{cor:reg_total_err_bound_matern}
	Let the assumptions of \Cref{thm:total_error_bound} hold and let $\lambda \geq 0$. Then, for the RKHS $\bH = H^{\nu + n/2}(\Omega)$ induced by the Matérn kernel $\kernel_{\nu}$, we have the pointwise error bound
	\begin{align}\label{eq:reg_total_err_matern}
		\| (\calK - \wh \calK_\lambda^\mathrm{MC}) f \|_\infty \leq C \left( h_\calX^{l+1/2} + \lambda^{1/2} \right) + \veps\cdot \sqrt{\frac{d}{\lambda}}
	\end{align}
	for every $\veps > 0$ and all $f \in H^{\sigma_{n,l}}$ with $\|f\|_{\bH} \leq 1$ with probability at least $1 - 2 d\exp\left( -\frac{m\varepsilon^2}{2}\right)$ and constant $C > 0$ only depending on $\kernel_{\nu}$.
\end{corollary}
\begin{proof}
	The proof is analogous to the second part of \Cref{cor:reg_total_err_bound} and leverages the uniform bound $\|\kernel_\nu\|_\infty = 1$. 
\end{proof}

\section{Balancing deterministic and probabilistic error}\label{app:numerical_considerations}
Here, we discuss the practical choice of the parameters $h_{\mathcal X}$, 
$m_{\mathrm{train}}$ and~$\lambda$ for the stochastic kernel EDMD method 
when using the Wendland $C^2$-kernel ($l=1$) in spatial dimension $n = 1$. The derivation for the Matérn kernel can be done in a similar fashion. 

The deterministic error associated with the projection onto the span of the kernel
translates, for Wendland smoothness $l=1$, into the estimate
\[
\| \mathcal K f - \widehat{\mathcal K} f \|_\infty
\le 
C
\|f\|_\bH
\, h_{\mathcal X}^{\,l + 1/2}
=
C
\|f\|_\bH
\, h_{\mathcal X}^{3/2}.
\]
from \Cref{cor:det_err_wendland}. Thus, decreasing fill distance improves the deterministic accuracy at rate $O(h_{\mathcal X}^{3/2})$. As $h_{\mathcal X}$ decreases, the number $d$ of design points increases, and the smallest eigenvalue of the kernel matrix $\mathbf K_{\mathcal X}$ decreases.
For the unregularized inverse, we start from the basic inequality $\max_{v\in\one} v^\top \mathbf K_{\mathcal X}^{-1} v \leq \frac{d}{\lambda_\mathrm{min}(\mathbf{K}_\calX)}$. A classical result from kernel interpolation theory provides a lower bound on the smallest eigenvalue via the \emph{separation distance} $q_\calX = \min_{i \neq j} \|x_i - x_j \|_2$ of the set $\calX = \{x_1, \dots, x_d\}$: Specifically, we have $\lambda_{\min}(\bbK_\calX) \geq C q_\calX^{2 s - n}$ with smoothness $s > n/2$ and dimension $n$ of the underlying kernel, see~\cite{Wenz2026, Scha1995, Wend04}. For Wendland kernels $\kernel_{n,l}$, where $s = \frac{n+1}{2} + l$, this gives rate $q_\calX^{2l+1}$. We assume equispaced points in dimension $n = 1$, for which $O(q_\calX) = O(h_\calX)$ holds, such that $\lambda_{\min}(\bbK_\calX) = O(q_\calX^3) = O(h_\calX^3) = O(d^{-3})$, and
\begin{align*}
	\left(
	\max_{v\in\one} v^\top \mathbf K_{\mathcal X}^{-1} v
	\right)^{1/2} = O\left( \sqrt{\frac{d}{d^{-3}}}\right) = O(d^2) = O(h_\calX^{-2}).
\end{align*}
Under Tikhonov regularization, \Cref{cor:reg_total_err_bound} shows that $\left(
\max_{v\in\one} v^\top \mathbf K_{\mathcal X,\lambda}^{-1} v
\right)^{1/2}
= O\!\left(\sqrt{\frac{d}{\lambda}}\right)$. However, regularization introduces an additional deterministic bias:
\[
\| \mathcal K f - \widehat{\mathcal K}_\lambda f \|_\infty
\le
C
\left(h_{\mathcal X}^{3/2} + \lambda^{1/2}\right)
\|f\|_{\bH}.
\]
To prevent the regularization bias from dominating, one requires $\lambda = O(h_{\mathcal X}^{3})$. 

Moreover, Monte Carlo estimation of $\mathbf K_{\calX, \calY}$ introduces a stochastic error. From \Cref{thm:mc_error1}, with probability at least
$1 - 2d e^{-m \varepsilon^2/(2\|\kernel\|_\infty)}$,
\[
\|\widehat{\mathcal K}f - \widehat{\mathcal K}^{\mathrm{MC}} f\|_\infty
\le
\varepsilon \,
C_{\mathcal X}
\|f\|_{\bH},
\]
where $C_\calX$ depends on the data points $\calX$. Solving the probability expression for~$\varepsilon$ gives $\varepsilon = O(m^{-1/2})$, again ignoring logarithmic terms in~$d$.
Since equispaced grids satisfy $d=O(h_{\mathcal X}^{-1})$, the Tikhonov control factor becomes
\[
\sqrt{\max_{v} v^\top \mathbf K_{\mathcal X,\lambda}^{-1} v}
= 
O\!\left(\sqrt{\frac{d}{\lambda}}\right)
= 
O\!\left(\frac{h_{\mathcal X}^{-1/2}}{h_{\mathcal X}^{3/2}}\right)
=
O(h_{\mathcal X}^{-2}).
\]
To ensure the probabilistic error is of the same order as the
deterministic error $O(h_{\mathcal X}^{3/2})$, we solve $\varepsilon \cdot O(h_{\mathcal X}^{-2}) = O(h_{\mathcal X}^{3/2})$, which gives $\varepsilon = O(h_{\mathcal X}^{7/2})$. Since $\varepsilon = O(m^{-1/2})$, this yields
\[
m_{\mathrm{train}}(h_{\mathcal X})
= O\!\left(h_{\mathcal X}^{-7}\right) = O(d^7).
\]
Since this is infeasible in practice, we conducted numerical experiments to test more moderate orders of growth of $m_\mathrm{train}$. We found the orders 
\begin{align*}
	m_{\mathrm{train}}(h_{\mathcal X})
	= O\!\left(h_{\mathcal X}^{-3/2}\right), \quad
	\lambda(h_{\mathcal X})
	= O\!\left(h_{\mathcal X}^{-3}\right),
\end{align*}
to produce satisfying balancing of the two error sources, as displayed in \Cref{fig:linf_ou_const_vs_adapt_left_right}.

\end{document}